\documentclass[12pt]{article}
\usepackage{amssymb,amsmath}
\usepackage{locenglish,math}
\usepackage{graphicx,pst-plot,psfig,pstricks}
\newcommand{\LastUpdate}{10 December 2001}
\newcommand{\Vo}{\oc{\BbV}}
\newcommand{\NTF}{\textsf{NonTrivFix(T)}}
\begin{document}
\begin{center}
\textbf{\large{\textsf{Random walks in random environment on trees\\
and multiplicative chaos}}}

\vskip1cm
\parbox[t]{14cm}{
Mikhail {\sc  Menshikov}$^a$ and Dimitri {\sc Petritis}$^b$\\
\vskip5mm   
{\scriptsize
\baselineskip=5mm
a. Department of Mathematical Sciences, University of Durham\\
South Road, Durham DH1 3LE, United Kingdom, Mikhail.Menshikov@durham.ac.uk\\
\vskip5mm
b. Institut de Recherche
Math\'ematique, Universit\'e de Rennes I and CNRS UMR 6625\\
Campus de Beaulieu, 35042 Rennes Cedex, France, Dimitri.Petritis@univ-rennes1.fr}}
\end{center}

\vskip1cm
{\small
\centerline{\LastUpdate}
\vskip1cm
\baselineskip=5mm
\noindent
{\bf Abstract:}
We study  random walks in a random environment on a regular, rooted, coloured tree. The asymptotic
behaviour of the walks is classified for ergodicity/transience in terms of
the geometric properties of the matrix describing the random environment. 
A related problem, with only one type of vertices and quite stringent conditions on the
transition probabilities but on general trees has   
been considered previously in the literature \cite{LyoPem}. In the presentation we give here,
we restrict the study of the process on a regular
graph instead of the irregular graph used in \cite{LyoPem}.
The close connection between various problems on random walks in
random environment and the so called multiplicative chaos martingale is underlined
by showing that the classification of the random walk problem can be drawn by the
corresponding classification for the multiplicative chaos, at  least for those
situations where both problems have been solved by independent methods. 
The chaos counterpart of the problem we considered here has not yet been solved.
The results we obtain for the random walk problem localise the position  of the
critical point. We conjecture that the additional conditions needed for the chaos
problem to have non trivial solutions will be the same as the ones needed
for the random walk to be null recurrent.\\
\vskip4mm
}

\section{Introduction}
\subsection{Notations}

Let $d$ be a fixed non-negative
integer. We consider the rooted
regular tree of order $d$, \textit{i.e.}\  a
connected graph without loops
with a denumerable set of
vertices $\mathbb{V}$ and a
denumerable set of non oriented
edges $\mathbb{A}(\mathbb{V})$. There is a distinguished
vertex called the root that has
degree $d$; all other vertices
have degree $d+1$. Vertices are
completely determined by giving
their genealogical history from
their common ancestor, the root; 
hence they are bijectively indexed by the
set of sequences of arbitrary length
over an alphabet of $d$ letters. We use
the same symbol for the indexing
set so that
$\mathbb{V} =\cup_{n=0}^{\infty} \mathbb{V}_n$ with 
$\mathbb{V}_0 = \{ \emptyset  \}$ and 
$\mathbb{V}_n = \{v = (v_1,\ldots,v_n) : v_i \in \{1,\ldots,d\}, i =
1\ldots n\}$ for $n \geq 1$. For every $v \in
\mathbb{V}$, we denote $|  v | $
the length of the path from $v$ to the
root \textit{i.e.} \ the number of edges
encountered. For $v \in
\mathbb{V}$ and $k \leq |  v | $
we denote by $v| _k$ the truncation
of the sequence $v$ to its $k$ first
elements, \textit{i.e.}  if $v = (v_1 \ldots  v_n)
\in \mathbb{V}_n$ and $k \leq
n$, then $v| _k = (v_1 \ldots  v_k)
\in \mathbb{V}_k$; the symbol
$v| _k$ must not be confused
therefore with $v_k$, representing
the letter appearing at the $k$-th
position of the sequence. For $0
\leq k < \ell \leq |v| $ we
denote $ v\vert_\ell^k $ the subsequence of length
$\ell-k$ defined by $ v\vert^k_\ell =
(v_{k+1},\ldots ,v_\ell)$. If $u \in
\mathbb{V}$, we write $u \leq v$ if
$|  u |  \leq |  v | $ and $v
= (u_1,\ldots  u_{|  u | },
v_{|  u| +1},\ldots  v_{|  v
| })$ \textit{i.e.}  if $u$ is the initial
sequence of $v$; we write $u<v$ when $u\leq v$ and $|u|<|v|$.
Similarly for
every  sequence $u$ and any
letter $\ell \in \{1\ldots d\}$, the
sequence $u\ell$ will have length
$|  u |  + 1$ and last letter
$\ell$.

Edges are unordered pairs $\langle u, v\rangle$
of adjacent vertices $u$ and $v$. Now
if $u$ and $v$ are the end vertices
of an edge, either $u \leq v$ with
$|  u |  + 1 = |  v | $ or $v
\leq u$ with $|  v |  + 1 = | 
u | $. In both cases, there is a
vertex, uniquely defined, that is the
most remote from the root among the
two end vertices of the edge. Since every
vertex has an unique ancestor, every
edge is uniquely defined by its
most remote vertex. Hence, every
vertex $v \in \mathbb{V} \backslash
\{ \emptyset  \}$ defines an edge 
$a(v) = \bra v\vert_{\vert\, v\,\vert-1},v\ket$. 
Edges are thus also indexed
by the set $\mathbb{V}$, more precisely by $\Vo = \mathbb{V} 
\setminus
\lbrace\emptyset \rbrace$ and we denote $a(v)$ 
the edge defined by $v$; therefore
$\mathbb{A}(\mathbb{V})\simeq\oc{\BbV}$. Since the set of vertices
$\mathbb{V}$ uniquely determines the set of edges, we use by abuse of
notation the symbol $\mathbb{V}$ to denote also the tree.

If $u, v  \in \BbV $ and $u < v$, we denote $[u; v]$
the (unique) path from $u$ to $v$
 \textit{i.e.}\  the collection of edges $(a_1, a_2, \ldots )$ with $a_j
\equiv a(v|_{|u|+j})$, for $j = 1, \ldots , |v| $. For every $u\in\BbV$, 
the symbol $[u;u]$ denotes an  empty set of edges. If $u$ and $v$ are not comparable
vertices, \textit{i.e.}\ neither $u\leq v$ nor $v\leq u$ holds, although there is
a canonical way to define the path $[u;v]$, this definition is not necessary in the
present paper and hence omitted.
We write simply $[v]$ to denote the path joining the root to $v$, namely $[\emptyset;v]$.

At every edge $a$ we assign a number $\xi_a \in [0, \infty[$
in some specific manner. This specification differs from
model to model and since various models are considered here,
we don't wish to be more explicit about these variables at
the present level. Mind however that 
the numbers
$(\xi_a)_{a\in\mathbb{A}(\mathbb{V})}$ are random
variables neither necessarily independent nor necessarily
equi-distributed. For the time being, we only assume that we
dispose of a specific collection
$(\xi_a)_{a\in\mathbb{A}}$, called the \textit{edge-environment}. 

\subsection{Multiplicative chaos}
Let $(\mathbb{V}, (\xi_a)_{a\in\mathbb{A}(\mathbb{V})})$ be
a given tree and a given edge environment.
For $u, v \in \mathbb{V}$, with $u < v$ we denote 
\[\xi[u;v]
= \prod\limits_{a\in[u;v]}\xi_a\]
the product of environment values encountered on the path of
edges from $u$ to $v$; the symbol $\xi[v]$ is defined to
mean $\xi[\emptyset ; v]$ and $\xi[v;v]$ --- as a product over
an empty set --- is consistently defined to be 1. It is not necessary for the
purpose of the
present article to define  the value of $\xi[u;v]$ when $u$ and $v$ are not comparable.

For every $u \in \mathbb{V}$,  we
consider the process $Y_n(u)_{n\in\mathbb{N}}$ defined by
$Y_0(u) = 1$ and 
\[Y_n(u) =
\sum\limits_{v\in\mathbb{V}_{n+|u|}; v>u}\xi[u;v]=
\sum\limits_{v\in\mathbb{V}_{n+|u|}; v>u}\prod\limits_{a\in[u;v]}\xi_a,\]
 for $n
\geq 1$. This process is known as the \textit{
multiplicative chaos process}.
Notice that even when $(\xi_a)_{a\in\mathbb{A}}$ is a
family of independent random variables, the random
variables $\xi[u;v]$ are not independent for $v$ scanning
the set $\mathbb{V}_{n+|  u| }$. 
Hence the asymptotic behaviour of $Y_n(u)$ when $n
\r  \infty$ is far from trivial and it is
studied for several particular cases of dependences of the
family $(\xi_a)$ in an extensive literature; see for
instance 
\cite{KahPey,DurLig,Gui,ColKou,Liu97,Liu98,WayWil}.

The study of the asymptotics of the process $(Y_n)$  
is done by various techniques. One such technique is
 by remarking that for $n \geq 1$, the process
can be written as
\begin{eqnarray*}
Y_n(u)
&=&\sum\limits_{v\in\mathbb{V}_{n+|u|}; v>u}\xi[u;v] 
= \sum\limits_{v\in\mathbb{V}_{n+|u|};v>u}
\xi[u;v|_{|u|+1}]\,\xi[v|_{|u|+1};v]  \\
&=&\sum\limits_{w\in\mathbb{V}_{|u|+1};w>u}
\xi[u;w]\,\sum\limits_{v\in\mathbb{V}_{|u|+n}; v>u}\,\xi[w;v]\\
&=&\sum\limits_{w\in\mathbb{V}_{|u|+1}; w>u}\xi[u;w]\,Y_{n-1}(w).
\end{eqnarray*}
If the limit $\lim\limits_{n\r \infty} Y_n(u) \ \elaw\  
Y(u)$ exists in distribution for all $u \in \mathbb{V}$ then
it must verify the functional equation
\begin{equation}
\label{eq:fun-eq}
Y(u)\  \elaw  
\sum\limits_{w\in\mathbb{V}_{|u| +1};w>u}\,\xi[u;w]\, Y(w).
\end{equation}
The process $(Y_n(u))_n$ and the corresponding functional equation
(\ref{eq:fun-eq})
are thoroughly studied in the literature for some particular
choices of dependencies of the family $(\xi_a)$. 

A second technique of study of the asymptotics is by martingale
analysis. We have, as a matter of fact,
\begin{eqnarray*}
Y_n(u)&=& \sum\limits_{v\in\mathbb{V}_{n+|u|}; v>u}\xi[u;v]\\
&=& \sum\limits_{v\in\mathbb{V}_{n+|u|}; v>u} \prod\limits_{a\in [u;v]} \xi_a\\
&=& \sum\limits_{v_1,\ldots,v_n\in \{1,\ldots,d\}}
\xi_{a(uv_1)}\xi_{a(uv_1v_2)}\ldots
\xi_{a(uv_1v_2\ldots v_n)}
\end{eqnarray*}
Now, if for any fixed $u\in\BbV$,  $(\cF_n^{(u)})$ denotes the natural filtration
$\cF_k^{(u)}=\sigma(\xi_{a(uv_1\ldots v_k)}, v_i\in\{1,\ldots,d\}, i=1,\ldots,k\}$ for
$k\in\BbN$, we have
\[\BbE(Y_n(u)|\cF_{n-1}^{(u)})= \sum\limits_{v_1,\ldots,v_n\in \{1,\ldots,d\}}
\xi_{a(uv_1)}\ldots\xi_{a(uv_1\ldots v_{n-1})}\BbE(\xi_{a(uv_1v_2\ldots v_n)}|\cF_{n-1}^{(u)})\]
and in the special case where the distribution of $\xi_{a(uv_1v_2\ldots v_n)}$ depends
solely on $v_n$ and the random variables are independent for different generations,
 the previous formula simplifies into
\[\BbE(Y_n(u)|\cF_{n-1}^{(u)})=Y_{n-1}(u) \sum_{v_n=1}^d \BbE(\xi_{a(uv_1v_2\ldots v_n)}).\]
In this special case and provided that $\sum_{v_n=1}^d \BbE(\xi_{a(uv_1v_2\ldots v_n)})$ is less, equal,
or more than 1, the process is a non-negative supermartingale,  martingale or
submartingale. It will be shown in section 
\ref{sect:mult-chaos-results} that it is enough
to consider the martingale case since by suitable renormalisation of the
random variables $(\xi_a)$ and $(Y_n)$ we can always limit ourselves in the
study of martingales.
 
Although the process $(Y_n)$ is thoroughly studied, 
the closely related process 
\[Z_n(u) = \sum_{k=0}^n
Y_k(u) \ \ \textrm{for}\ \  n \geq 0\]
 does not seem --- to the best of our
knowledge --- to have attracted much attention. However, if
we are interested in connections between multiplicative
chaos and random walks in random environment on a tree, it
is this latter process that naturally appears in both
subjects.

The relevant question that can be addressed for the
multiplicative chaos problem is whether the functional
equation (\ref{eq:fun-eq}) has a non trivial solution.
The precise statement of this question is 
model-dependent and some instances of it will be examined below, in subsection
\ref{models}.

\subsection{Nearest neighbours random walk on a tree in an inhomogeneous
environment}
To every vertex $u=(u_1,\ldots,u_{|u|}) \in \Vo$ are assigned $d+1$
numbers $(p_{u,0}, p_{u,1}, \ldots  , p_{u,d})$ with
$p_{u,0}>0$, 
$p_{u,i} \geq 0$ $\forall i = 1, \ldots,  d$ and $\sum_{i=0}^d$
$p_{u,i} = 1$. To $u \in \mathbb{V}_0=\{\emptyset\}$ are assigned only $d$
numbers $(p_{\emptyset,1}, \ldots  p_{\emptyset,d})$ with $p_{\emptyset,i} \geq 0$
$\forall i = 1, \ldots,  d$ and $\sum_{i=1}^d$ $p_{\emptyset,1} = 1$.
These numbers will be random variables with some
specific dependence properties that will be defined later.
These numbers stand for transition probabilities of a
reversible Markov chain $(X_n)_{n\in\mathbb{N}}$ on the tree
verifying for $|u|  \geq 1$ 
\[P_{u,v} = \mathbb{P}(X_{n+1} =v|X_n = u)=\left\{
\begin{array}{ll}
p_{u,0} & \textrm{if}\ \ v=u|_{|u|-1}\\
p_{u,v_{|v|}} & \textrm{if}\ \ u=v|_{|v|-1}\\
0 & \textrm{otherwise.}
\end{array}\right. \]
For $u = (\emptyset )$ we have the slightly modified
transition probabilities  
\[P_{\emptyset ,v} = \mathbb{P}(X_{n+1}
= v |  X_n = (\emptyset )) =
\left\{
\begin{array}{ll}
p_{\emptyset, v_1} & \textrm{if}\ \ v\in\mathbb{V}_1\\
0 & \textrm{othrewise.}
\end{array}\right. \]

For $u \in \mathbb{V}$ with 
$|u|  \geq 2$ we consider the edge 
$a(u) = \bra u| _{|u| -1}, u\ket$ 
and attach to this edge the variable
\[\xi_{a(u)} = {\frac{p_{u|_{|u|-1},u_{|u|}}}
{p_{u|_{|u|-1},0}}} \in [0,\infty[.\]
For $u \in \mathbb{V}_1$ we attach $\xi_{a(u)} =
p_{\emptyset ,u_1}$.
Notice that we have not required irreducibility for the
Markov chain $(X_n)_{n\in\mathbb{N}}$. The condition $0 \leq
\xi_{a(u)} < \infty$ only excludes vanishing of  the probability of
transition towards the root on every active branch
of the tree. 
We require also some mild integrability conditions
on the random variables in the form $\BbE \xi_{a(v)}\log^+(\xi_{a(v)})<\infty$, where
for all $z>0$, $\log^+(z)=\max(0, \log(z))$.
One can easily check the validity of the following
\begin{lemm}\label{lem:inv-meas}  
For every $v \in \BbV$ define the variable
\[\pi[v] = \left\{\begin{array}{ll}
\pi[\emptyset ]\, \xi[v] \frac{1}{p_{v,0}}&
\textrm{if}\ \ v\in\Vo\\
\pi[\emptyset] &\textrm{if}\ \ v=(\emptyset),\end{array}\right.
\] 
with
$\pi[\emptyset ]$ an arbitrary constant. Then $\pi[v]$ verifies the
stationarity condition
\[\sum\limits_{v\in\mathbb{V}} \pi[v] P_{v,v'} = \pi[v'], \quad\forall
v'
\in \BbV.\] 
\end{lemm}
We assume that the variable $(p_{v,0})^{-1}$ is well behaved
uniformly in $v$. To avoid technical difficulties, we assume that
\begin{equation}
\label{inverse-integrability}
\BbE((p_{v,0})^{-1})<\infty.
\end{equation}
Then, apart  the factor
$\frac{1}{p_{v,0}}$,  the expression for the invariant measure
$\pi[v]$ involves the product $\xi[v]$ of variables along
the edges of the path form $\emptyset $ to $v$ as  was the case
in the expression of multiplicative chaos. The form of the
invariant measure established in lemma \ref{lem:inv-meas} has been
already established in \cite{LyoPem}, where the
problem of random walk in a random environment on an
inhomogeneous tree was studied in a particular case of
dependence of the random variables $(\xi_a)$. 
However, the
close analogy with multiplicative chaos,
although reminiscent of considerations in 
\cite{Lyo92}, does not seem to have been
exploited.

The relevant questions that can be addressed in the context
of random walk in inhomogeneous environments concern ergodicity or
non-ergodicity, recurrence or transience that arise for
almost all realisations of the variables $(\xi_a)$.

\Rk Notice that since our chains are periodic of period 2, we use the term \textit{ergodic}
to mean positive recurrent chains (or else admitting a stationary probability distribution $\pi$)
and not in the sense that $\BbP(X_n=v|X_0=u)$ converges towards $\pi(v)$; as a matter of fact this
conditional probability does not admit a limit when $n\rightarrow \infty$, it is only the subsequence 
of even  (resp.\ odd) times that converges towards this limit according to the
parity of the difference $|v|-|u|$, terms of complementary parity being always 0.

\subsection{Models covered by the present formalism}
\label{models}
We present below a unified treatment of 
both the multiplicative chaos
process and the random walk problem stating in the same theorem the asymptotic
behaviour of the limiting chaos
process and of the random walk. Several models fit the present formalism; by making appropriate identifications
of random variables, the random walk in random environment on $\mathbb{N}$ or the 
problem of random strings in a random environment can be rephrased in the present language.

\subsubsection{Random walk in a random environment on a regular tree}
At every vertex $v\in\Vo$ is assigned a $(d+1)$-dimensional random vector with
positive components $(p_{v,0},\ldots,p_{v,d})$ verifying $\sum_{j=0}^dp_{v,j}=1$.
For the vertex $v=\emptyset$, the corresponding random vector is $d$-dimensional
and its components verify $\sum_{j=1}^d p_{\emptyset,j}=1$. These random vectors
are independent for different $v$'s and, for $v\in\Vo$ they have the same distribution.
The random walk in the random environment defined by these random vectors on the regular
tree is the Markov chain with transition matrix
\[\BbP(X_{n+1}=v|X_n=u)=\left\{\begin{array}{ll}
p_{u,j} & \textrm{if}\ \ u\in\BbV \ \ \textrm{and}\ \ v=uj, j=1,\ldots,d\\
p_{u,0} & \textrm{if}\ \ u\in\Vo \ \ \textrm{and}\ \ v=u|_{|u|-1}\\
0 & \textrm{otherwise.}
\end{array}\right.\]
Notice again that although the random vectors $(p_{v,0},\ldots,p_{v,d})$ for $v\in\Vo$ are independent,
their components \textit{cannot} be independent since they satisfy $\sum_{j=0}^dp_{v,j}=1$.
Assigning to every edge the ratio of outwards over inwards probabilities, the problem can
now be rephrased in a way that it fits the general formalism.
Let $\bom{\eta} = (\eta_1 , \ldots , \eta_d)$,
 be a vector of non-negative
random variables $\eta_i$, $i=1,\ldots,d$, 
having the same distribution with $p_{v,i}/p_{v,0}$, for $v\in\Vo$, 
with not necessarily independent nor identically
distributed components. We assume the law of the random vector is
explicitly known with 
$\mathbb{E} \eta_i < \infty$ and 
$\mathbb{E} \eta_i \log^+\eta_i < \infty,\  
\forall i = 1,
\ldots, d.$  Moreover, to avoid technicalities we assume that although
the support of the random variables $\eta_i$ extends up to 0, their
law has no atom at 0.

To the edge $a(v)$, having most remote vertex $v \in \Vo$, we assign
the random variable $\xi_{a(v)}$ having the same distribution as
$\eta_{v_{|v|}}$; the 
variables $\xi_{a(v)}$ and $\xi_{a(v')}$ are
independent if $v|_{|v|-1} \not = v'|_{|v'|-1}$. Notice that if the
components of the random vector $\eta$ are not independent, the
variables $\xi_{a(v)}$ and $\xi_{a(v')}$ with $|v| = |v'|$ and
$v|_{|v|-1} = v'|_{|v'|-1}$ are not independent either.

When $\bom{\eta}$ has independent and identically distributed
components, the process $(X_n)$ has been studied in the context of
multiplicative chaos and the existence of non trivial solutions of
the functional equation (\ref{eq:fun-eq}) is established in
\cite{KahPey,Gui,DurLig}. In a context of statistical mechanics this
functional equation is also studied in \cite{ColKou,Pet} and in a
more general branching process in \cite{GraMauWil,Liu97}. For $\bom{\eta}$
with components having a general joint distribution, the process
$(Y_n)$ has been studied in \cite{DurLig} and when the dimension
$d$ of the vector is also a random variable having a general joint
distribution with the components of $\eta$ has been studied in
\cite{Liu98} (see \cite{LiuHAB} for a recent survey of the state of
the art on the subject).

The results are expressed in terms of the functions 
\[f(x) =
\mathbb{E}\;(\sum_{i=1}^d \eta_i^x), x \in \mathbb{R}^+
\ \ \textrm{and} \ \  g(x)=\log f(x),\] 
and of the parameter $\lambda = \inf\limits_{x\in[0,1]} f(x).$ 

The equation (\ref{eq:fun-eq}) reduces asymptotically into 
the  form
\begin{equation} \label{eq:fun-eq-rwre}
Y[\emptyset] \elaw \sum\limits_{v\in\mathbb{V}_1}\; \eta_{v_1}
Y'[v]
\end{equation}
where $(Y'[v])_{v\in\mathbb{V}_1}$ are mutually independent
variables and independent of $(\eta_i)$ each having the same
distribution as $Y[\emptyset] =
\lim\limits_{n\rightarrow\infty} Y_n[\emptyset].$
Since only equality in distribution is required in equation
(\ref{eq:fun-eq-rwre}), the variables $Y[\emptyset]$ and $Y'[v]$
can be chosen independent.

The random walk problem in a random environment on a tree has
been first considered in several papers by 
Lyons, Pemantle, and Peres 
(see \cite{Lyo92,LyoPem,PemPer} for instance.) 
In particular in \cite{LyoPem}
Lyons and Pemantle studied the case of more general trees than the
ones considered here, namely trees whose degree is not
constant but have merely a finite  branching
number. On the other hand, these authors make more stringent
assumptions on the distribution of the random variables. 
We provide here proofs based on results on 
multiplicative chaos, 
totally independent from the proofs 
given by these authors.
We
state here our result for the case 
of a regular tree (constant
degree) but with \textit{arbitrary distribution} 
for the components of
the random vector $\eta_i$ \textit{i.e.}\ with  
neither independent nor
identically distributed components. It is worth noticing that
although we limit ourselves in the case of regular trees, this limitation
is not so crucial. As a matter of fact, the recent results on multiplicative chaos
allow the treatment of random branching numbers and not to excessively burden
the present paper, this will be considered in a subsequent 
article.
\begin{theo}
\label{th-Lyo-Pem}
Let $\lambda = \inf\limits_{x\in[0,1]} f(x)$ and
$x_0 \in [0,1]$ be such that $f(x_0) = \lambda.$
Then
\begin{enumerate}
\item If $\lambda < 1$, then almost surely the random walk is ergodic and
$Z_\infty < \infty$.
\item If $\lambda > 1$, then almost surely the random walk is transient,
$Y_\infty=\infty$, and $Z_\infty = \infty$. 
\item If $\lambda = 1$ and moreover $f'(1) < 0$, 
then almost surely 
$0<Y_\infty<\infty$, $Z_\infty=\infty$, and the random walk is null-recurrent.
\end{enumerate}   
\end{theo}

Assertions 1.\ and 2.\ in this theorem can be viewed as special cases of the
corresponding assertions of theorem \ref{main-th}. However, we shall prove this theorem
independently, using a new method, because we wish to demonstrate how
techniques developed in the context of multiplicative chaos can be used in random
walk problems. Moreover, this result is interesting in its own because it is
a generalisation of the result in \cite{LyoPem}, valid for random variables
that need not be either independent or identically distributed.
Although this result is presented as a single theorem, it follows from
several partial theorems on 
multiplicative chaos 
stated and/or proved in section \ref{sect:mult-chaos-results}.

\subsubsection{Random walk in a random environment on a coloured tree}
This problem is reminiscent of the problem on random strings in a random environement,
studied in
\cite{ComMenPop}, where non reversible Markov chains on the tree $\mathbb{V}$
are considered and general conditions for transience/null
recurrence/ergodicity  are
given in terms of Lyapuonov exponent of a product
of matrices. Under some simplifying conditions the string problem can be rephrased
to fit 
the present formalism.
To describe the problem of random strings in a random environment, the underlying tree
must be enlarged to distinguish the $d$ children
of every vertex; this can be done by assigning  a colour index, chosen without replacement
from the set $\{1,\ldots,d\}$, to each child. For instance, a natural way to do so
is by assigning the colour $i\in\{1,\ldots,d\}$ to the vertex $v\in\Vo$ if
$v_{|v|}=i$.  The root is assigned an arbitrary colour
$\alpha\in\{1,\ldots,d\}$. Consequently, every edge $a(v)$ with $v\in\Vo$ is assigned
the bicolour $(ij)\in\{1,\ldots,d\}^2$, where $i=v_{|v|-1}$ and $j=v_{|v|}$.

In the sequel, we use the same colouring procedure to distinguish among children and the model
we consider is defined as follows.
To every $v\in\Vo$, assign the
$(d+1)$-dimensional vector with positive components $(p_{v,v_{|v|},j}, j\in\{0,\ldots,d\})$,
verifying $\sum_{j=0}^d p_{v,v_{|v|},j}=1$ for every $v\in\Vo$. These vectors are independent for
different $v\in \Vo$ but not equidistributed, 
their distribution depending on the colour $v_{|v|}\in\{1,\ldots,d\}$.
For the root $v=\emptyset$, we have the special assignment
$(p_{\emptyset, \alpha, j}, j\in \{1,\ldots,d\})$ with $\sum_{j=1}^d p_{\emptyset,\alpha,j}=1$.
The random walk in this random environment on the regurar tree is the
Markov chain with transition matrix
\[\BbP(X_{n+1}=v|X_n=u)=\left\{\begin{array}{ll}   
p_{u, u_{|u|}, j} & \textrm{if} \ \ u \in \BbV, v=uj, j\in\{1,\ldots,d\}\\
p_{u, u_{|u|}, 0} & \textrm{if} \ \ u \in \Vo, v=u|_{|u|-1}\\
0  & \textrm{otherwise.}
\end{array}\right.\]

Passing to the edge-indexed ratio of outwards over inwards probabilities, the model
can be rephrased to fit the present formalism.
   Let 
\[\boldsymbol{\eta} = \begin{pmatrix}
\eta_{11} & \cdots & \eta_{1d} \\
\vdots \\
\eta_{d1} & \hdots & \eta_{dd}
\end{pmatrix}\]
be a matrix of non-negative random elements of known joint distribution. 
The matrix
elements are not necessarily independent. 

The root is arbitrarily assigned a colour 
$\alpha \in \{1, \ldots,d \}$.
Every edge of the tree, indexed by its most remote vertex $v
\in \Vo$, carries thus a bicolour-type $(i,j)$ with $v_{|v|-1} =
i$ and $v_{|v|} = j$; it is assigned a
random variable  $\xi_{a(v)}$ that is a copy of
the random variable  $\eta_{ij}$. The random variables  $\xi_{a(v)}$, $\xi_{a(v')}$ 
are independent if
$v|_{|v|-1} \not =  v'|_{|v'|-1}$.
Thus, the random variables
attached to children of different parents (\textit{i.e.}\  non siblings) are
independent.

It is evident that the process $Y_n$ will depend on the colour
$\alpha$ attached to the root. To recall this, we denote
$Y_n^{(\alpha)}$ for the process starting from an
$\alpha$-coloured root.

Then 
\begin{eqnarray}
Y_n^{(\alpha)} [\emptyset]
&=& \sum\limits_{v\in\mathbb{V}_n} \xi_{a(v|_1)} \ldots
\xi_{a(v|_{n-1})}\xi_{a(v)}\\
&\elaw& \sum\limits_{v\in\mathbb{V}_1} \eta_{\alpha v_1}
Y_{n-1}^{(v_1)} [v],
\end{eqnarray}
where $Y_{n-1}^{(v_1)} [v]$ are mutually independent random variables for
different $v_1$ and independent of the random variables $\eta_{\alpha v_1}$. 
If the sequence $(Y_n^{(\alpha)} [\emptyset])_n$ has a 
limit (in probability)
when $n\r\infty$, then $\lim_{n\r\infty}Y_n^{(\alpha)} [\emptyset]$
will depend only on $\alpha$; similarly, for every fixed value of $\alpha$, 
the $\lim_{n\r\infty}Y_{n-1}^{(v_1)} [v]$ will depend only on $v_1$. 
Therefore, if the limit $Y_\infty^{(\alpha)}(u)$ exists $\forall
\alpha \in \{1,\ldots,d\}$ and $\forall u \in \mathbb{V}$ 
it must verify
the functional equation
\begin{equation}
\label{eq:fun-eq-col-tree}
Y^{(\alpha)} \ 
\elaw \ 
\sum_{\beta=1}^d
\eta_{\alpha\beta} Y^{(\alpha,\beta)}\quad \forall \alpha\in\{1,\ldots,d\},
\end{equation}
where $Y^{(\alpha,\beta)}$ are mutually independent random variables
and independent of the random variables $\eta_{\alpha\beta}$. For every
$\alpha,\beta \in\{1,\ldots,d\}$ the random variables $Y^{(\alpha,\beta)}$
and $Y^{(\beta)}$ have the same law; moreover, since equation (\ref{eq:fun-eq-col-tree})
holds only in distribution, they can be chosen independent.
A very particular case of the equation (\ref{eq:fun-eq-col-tree}) has been
considered in \cite{BenNas}.
To the extend of our knowledge, the functional equation (\ref{eq:fun-eq-col-tree})
in
this generality has not
been studied in the literature. Restricting ourselves in the case of
integrable solutions and taking expectations from both sides,
it is evident that a necessary condition for the existence of
non-trivial solution is that $1$ is an eigenvalue of the matrix
$\mathbb{E}\eta$. This condition seems far from being sufficient.
Some conjectures about this functional equation 
are formulated in the
last section of this paper.

Let's consider the random walk in random environment  problem. Our main result is formulated in
the following
\begin{theo}\label{main-th}
Let $\bom{m}(x) = 
\begin{pmatrix}
\BbE(\eta_{11}^x) & \cdots & \BbE(\eta_{1d}^x) \\
\vdots \\
\BbE(\eta_{d1}^x) & \hdots & \BbE(\eta_{dd}^x)
\end{pmatrix}$      for $x \in [0,1]$.
Assume that the matrix $\bom{m}(x)$ 
is regular \textit{i.e.}\
there exists some
integer $N$ such that  that for every $x\in[0,1]$, 
$(\bom{m}(x)^N)_{ij} > 0\quad \forall i,j$. 
Denote
by $\rho(x)$ the largest eigenvalue of $\bom{m}(x)$ for $x \in [0,1]$ and
$\lambda = \inf\limits_{x\in[0,1]} \rho(x)$.
\begin{enumerate}
\item If $\lambda < 1$ the random walk is almost surely ergodic and $Z_\infty<\infty$ almost surely
\item If $\lambda > 1$ the random walk is almost surely transient and $Y_\infty=\infty$ almost surely.
\end{enumerate}
\end{theo}

The sections \ref{sect:ergodicity}, \ref{sect:technical}, and \ref{sect:proof}
below are devoted to the proof of the theorem \ref{main-th}.
In some places we used techniques 
already used in \cite{LyoPem},
either literally or in an extended version.
In other places, some independent 
probabilistic techniques are used that beyond their 
interest for the random walk 
problem they illuminate aspects
of the yet unsolved chaos problem.

\section{Ergodicity of the random walk in random environment for the case $\lambda <1$}
\label{sect:ergodicity}
The  proof of ergodicity is much simpler than the proof of non ergodicity.
\begin{lemm} 
\label{lem:ergod}
If $\rho(x) < 1$ for some 
$x \in [0,1]$ then $Z_\infty
= \sum\limits_{v\in\mathbb{V}}\, \xi[v] < \infty$  a.s.
\end{lemm}
\Proof
Fix first some $v \in \mathbb{V}_n$, $v = (v_1, \ldots, v_n)$ and
assign colour $\alpha$ at the root. Denote 
$\mathbb{E}_\alpha(\cdot)$
expectations for starting colour $\alpha$. We have
$\mathbb{E}_\alpha(\xi^x[v]) = m_{\alpha v_1}(x) m_{v_1v_2}(x)
\ldots m_{v_{n-1}v_n}(x)$.
Hence 
\[\mathbb{E}_\alpha(\sum\limits_{v\in\mathbb{V}_n} \xi^x[v]) =
\bom{e}_\alpha \bom{m}^n(x) 
\bom{e},\]
where the vectors $\bom{e}_\alpha$ and $\bom{e}$ have components $e_{\alpha i}=\delta_{\alpha i}$ 
and $e_=1$ respectively for $i=1, \ldots, d$.
If for some $x_0 \in [0,1]$, $\rho(x_0) < 1$ then
$\sum_{n=1}^\infty \bom{m}^n(x_0)$ is a well defined $d \times d$
matrix whose all elements are finite. Hence, for this
particular $x_0$, $\sum\limits_{v\in\Vo} 
\mathbb{E}_\alpha
\xi^{x_0}[v] < \infty$ and since $\xi^{x_0}[v] \geq 0$ this means
that $\sum\limits_{v} \xi^{x_0}[v] < \infty$ almost surely. 
The last majorisation means that there exists at most a
\textit{finite} subset $V \subseteq \mathbb{V}$ on which
$\xi^{x_0}[v] > 1$ for $v \in V$.
We have then:
$\sum\limits_{{v} \in \mathbb{V}}\xi[v] =
\sum\limits_{{v} \in V}\xi[v] + \sum\limits_{{v} \in
\mathbb{V}\setminus V} \xi[v]$ on $\mathbb{V}$.
On $\mathbb{V}\setminus V$, $ \xi^{x_0}[v]\geq \xi[v]$
 for $x_0\leq 1$ since $\xi[v]\leq 1$. Hence
\[ \sum\limits_{v\in\mathbb{V}}\xi [v]\leq \sum\limits_{v\in V}\xi[v]+
\sum\limits_{v\in\mathbb{V}\setminus V}\xi^{x_0}[v]\]
and the results follows immediately, both terms of the right hand side being finite.\eproof


\begin{coro}
If $\lambda <1$ the random walk in random environment  is almost surely ergodic.
\end{coro}

\Proof
It is enough to show that 
$\sum\limits_{v\in\BbV}\pi[v]<\infty$ since then, choosing the arbitrary normalising
value $\pi[\emptyset]$, we can guarantee that $\sum\limits_{v\in\BbV}\pi[v]=1$
showing that the random walk admits an invariant probability.
Now $\sum\limits_{v\in\BbV}\pi[v]=\pi[\emptyset]+\sum\limits_{n\in\BbN^*}\sum\limits_{v\in\BbV_n}\pi[v]$.
For $x_0\in]0,1[$ such that $\rho(x_0)<1$, we have
\begin{eqnarray*}
\BbE_\alpha(\sum\limits_{v\in\BbV}\pi[v]^{x_0})&=&
\BbE_\alpha(\pi[\emptyset]^{x_0})+\sum_{n\in\BbN^*}\sum\limits_{v\in\BbV_n}
\BbE_\alpha(\pi[\emptyset]^{x_0}\xi[v]^{x_0}p_{v,0}^{-x_0})\\
&=&\BbE_\alpha(\pi[\emptyset]^{x_0})+\sum_{n\in\BbN^*}\sum\limits_{v\in\BbV_n}
\BbE_\alpha(\pi[\emptyset]^{x_0}\xi[v]^{x_0})\BbE(p_{v,0}^{-x_0}),
\end{eqnarray*}
because $p_{v,0}$ is independent of $\xi[v]$ for all $v\in\oc{\BbV}$.
Now for $n\in\BbN^*$, $v\in\BbV_n$, and $x_0<1$, we have
$1\leq \frac{1}{p^{x_0}_{v,0}}\leq \frac{1}{p_{v,0}}$. Therefore,
$\BbE(\frac{1}{p^{x_0}_{v,0}})\leq \max_{v\in\BbV_n}\BbE(\frac{1}{p_{v,0}})
=\max\limits_{i=1,\ldots,d} \BbE(\frac{1}{p_{v'i,0}})$
for $v'=v|_{|v|-1}$ since the random variables 
$p_{v,0}$ and $p_{u,0}$ are independent and identically distributed
for $u,v\in\BbV_n$ with $u|_{|u|-1}\ne v|_{|v|-1}$.
Hence by using condition (\ref{inverse-integrability})
and the integrability condition
$\BbE_\alpha(\sum\limits_{v\in\BbV}\pi[v]^{x_0})<\infty$,
the proof is  completed by the same arguments as those used
in the proof of the previous lemma.
\eproof

\section{Technical results}
\label{sect:technical}
\subsection{Colouring}
Our results can be expressed more easily by
introducing a  systematic way to deal with
the colouring of the tree.

We call \textit{colouring function} 
a map $c:\BbV\r C$, where $C$ is a finite discrete set, 
\textit{the colour set}. We shall use 
in the sequel the colouring function with colour 
set $\{1,\ldots,d\}$ defined by 
\[c(u)=\left\{\begin{array}{ll}
u_{|u|} & \textrm{if}\ \ u\in \oc{\BbV}\\
\alpha\in\{1,\ldots,d\}  &  \textrm{if}\ \ u=(\emptyset).
\end{array}\right.\]
Thus the colouring of the tree introduces a natural ordering among vertices that are siblings;
 the root vertex that is assigned some arbitrarily prescribed colour
$\alpha$. To avoid ambiguity, we use --- when necessary --- subscripts
or superscripts to indicate the root colour, \textit{e.g.}\
$\BbE_\alpha(\cdot)$ for expectations, $Y_n^{(\alpha)}$ for
random variables, etc.

For $i,j\in\{1,\ldots,d\}$, we denote 
\[N_{ij}[u;v]=\sum_{k=|u|}^{|v|-1}\id_{(ij)}(c(v|_k),c(v|_{k+1}))\]
the number of edges in the path $[u;v]$ with
\textit{chromatic type} $(ij)$. If $u=(\emptyset)$, we
use as usual the simplified notation $N_{ij}[v]$.
Using the independence of the random variables encountered
on the edges of every path $[u;v]$, we remark that for every 
$x\in[0,1]$,
\begin{eqnarray*}
\BbE_\alpha(\xi^x[u;v])&=&m_{11}(x)^{N_{11}[u;v]}
\ldots 
m_{dd}(x)^{N_{dd}[u;v]}\\
&=& \prod_{i=1}^d \prod_{j=1}^d 
m_{ij}(x)^{N_{ij}[u;v]},
\end{eqnarray*}
where $m_{ij}(x)$, for $i,j=1,\ldots,d$ are the matrix elements of the matrix
$\bom{m}(x)$.
\begin{nota}{}
Let $\bom{\beta}=(\beta_{ij})_{i,j,\in\{1,\ldots,d\}}$ be
a $d\times d$ matrix with $\beta_{ij}\in [0,1]$ and 
$\sum_{i,j=1}^d\beta_{ij}=1$. Let $x\in[0,1]$. We introduce the
following functions
\begin{eqnarray*}
\psi(\bom{\beta})&=&\prod_{i=1}^d\prod_{k=1}^d
\left(\frac{\sum_{j=1}^d \beta_{ij}}{\beta_{ik}}\right)
^{\beta_{ik}},\\
\chi(x,\bom{\beta})&=&\prod_{i=1}^d\prod_{j=1}^d
m_{ij}(x)^{\beta_{ij}}, \textrm{and}\\
\phi(x,\bom{\beta})&=&\psi(\bom{\beta})\chi(x,\bom{\beta}),
\end{eqnarray*}
with the convention $0^0=1$.
\end{nota}

We are interested in the set of paths having pre-assigned chromatic type.

\begin{defi}{}
Let $\bom{n}=(n_{i,j})_{i,j\in\{1,\ldots,d\}}$ be a given
$d\times d$ matrix
of non-negative integer entries. Denote
$n=\sum_{i=1}^d\sum_{j=1}^d n_{ij}$ and define, on the $\alpha$-rooted tree
\[\BbV(\bom{n})=\{v\in\BbV_n: N_{ij}[v]=n_{ij}, \forall i,j=1,\ldots,d \}.\]
The matrix $\bom{n}$ is called \textit{admissible chromatic matrix} if
$\BbV(\bom{n})\ne\emptyset$.
\end{defi}

\Rk Notice that the matrix $(n_{ij})$ to be admissible must verify some obvious
constraints.  For a given path $[v]$, it is immediate that
\[\sum_{i,j}N_{ij}[v]=\sum_{k=0}^{|v|-1}\sum_{i,j}\id_{(ij)}(c(v|_k)c(v|_{k+1}))=|v|.\]
Hence the number 
\begin{equation}
\label{eq:constraint2} n=\sum_{i=1}^d\sum_{j=1}^d n_{ij}
\end{equation}
 uniquely defines the length of the path.
Similarly  
\begin{eqnarray*}
\sum_{i}N_{ij}[v]&=&\sum_{k=0}^{|v|-1}\sum_{i}\id_{(ij)}(c(v|_k)c(v|_{k+1}))\\
&=& \sum_{k=0}^{|v|-1} \id_{j}(c(v|_{k+1}))\\
&=&N_j[v],
\end{eqnarray*}
where $N_j[v]$ is the number of vertices of colour $j$ in the path $[v_1;v]$.
Hence, if $n_j=\sum_i n_{ij}$ and $k_j$ denotes the number of vertices of colour $j$ in the
path $[v]$, we have $k_j=n_j+\delta_{j, \alpha}$.
Now reading the chromatic sequence of the path $[v]$ in reverse order,
\textit{i.e.}\ starting from level $|v|$ down to the root, is equivalent in
transposing the matrix $\bom{n}$. Hence
$k_j=\sum_i n_{ji}+\delta_{j, v_n}$.
Therefore, the matrix $\bom{n}$ must verify the \textit{quasi-symmetry condition}
\begin{equation}
\label{eq:constraint1}
\forall i, \sum_{j=1}^d(n_{ij}-n_{ji})=\delta_{i,\alpha}- \delta_{i,v_n}
\end{equation}
to be admissible as chromatic matrix of a path $[v]$ in a $\alpha$-rooted tree.

\begin{lemm}
\label{lem:comb-estimates}
Let $\bom{n}=(n_{i,j})_{i,j\in\{1,\ldots,d\}}$ be a 
$d\times d$ matrix
of non-negative integer entries that is admissible as chromatic matrix of a 
$\alpha$-rooted tree. We denote $n_i=\sum_{j=1}^dn_{ij}$ and
$n=\sum_{i=1}^d n_i$.
Let $\BbV(\bom{n})=\{v\in\BbV_n:N_{ij}[v]=n_{ij}\}$
and $K(\bom{n})=\card \BbV(\bom{n})$. Then
\[K(\bom{n})=C_{n_1}^{(n_{11},\ldots,n_{1d})}\ldots
C_{n_d}^{(n_{d1},\ldots,n_{dd})},\]
where
\[C_{a_1+\ldots+a_d}^{(a_1,\ldots,a_d)}=
\frac{(a_1+\ldots+a_d)!}
{a_1!\ldots a_d!},\]
for non-negative integers $a_1,\ldots,a_d$
\end{lemm}
\Proof
We know that $\BbV(\bom{n})\ne \emptyset$ since $\bom{n}$ is admissible.
We must enumerate the paths $[v]$ with $v\in\BbV_n$ and having root colour
$\alpha$. The number $n_i$ represents the total number of edges in the path
with their first vertex of colour $i$. 
To describe a path $[v]$ with $N_{ij}[v]=n_{ij}$ for all
$i$ and $j$, it is enough to choose among the $n_i$ vertices
of colour $i$ their children of colour $1$ 
(there are $n_{i1}$ of them), \ldots, their children of colour $d$ (there are $n_{id}$
of them.) Hence for every colour $i$, there are 
$C_{n_i}^{(n_{i1},\ldots,n_{id})}$ choices and 
the choice for every $i$ is independent from the choice
for every $j\ne i$. \eproof

It is immediate to show by independence that 
\[\sum_{v\in\BbV(\bom{n})}\BbE_\alpha(\xi^x[v])=
K(\bom{n}) m_{11}(x)^{n_{11}}\ldots m_{dd}(x)^{n_{dd}}.\]

Notice that $\BbV(\bom{n})\subset V_n$ 
is a \textit{finite} set of vertices. 
We shall be interested into admissible chromatic matrices having large
$n$ and to the  related notion
of \textit{infinite} subtree $\BbU(\bom{n})\subset
\BbV$ constructed  by some iterative procedure to be described below.

First observe that if we denote $\bom{n}^{(\alpha)}$ an admissible chromatic
matrix for a path in a $\alpha$-rooted tree, the chromatic matrix
for the \textit{same} path in a $\alpha'$-rooted tree, call it $\bom{n}^{(\alpha')}$,
is such that
$\bom{n}^{(\alpha)}-\bom{n}^{(\alpha')}$ is a matrix having only one non-zero
column and this column has only 2 non-zero elements, these elements having
respectively values $\pm 1$ and $\mp 1$. Therefore, in the large
$n$ limit, when at least one element of  the $d\times d$ matrix $\bom{n}^{(\alpha)}$ is
necessarily of $\cO(n)$, modifying the root colour changes only two elements
of the matrix by $\cO(1)$. On the space $\textsf{ACM}_n$ of admissible chromatic matrices
with elements summing  to $n$, define the equivalence
relation $n\sim n'$ to mean $n_{ij}-n'_{ij}=\cO(1)$ for all $i,j=1,\ldots,d$.
Then obviously admissible chromatic matrices for any root colour are in the same
equivalence class.
Given a root colour $\alpha$, fix some $\bom{n}=\bom{n}^{(\alpha)}\in\textsf{ACM}_n$. 
Each $v\in  \BbV(\bom{n})$ defines a path $[v]$
from $(\emptyset)$ to $v$, so that vertices
$v\in  \BbV(\bom{n})$ can be viewed as leaves of a finite
tree. Now each such leaf, having colour $\alpha'=v_n$, can be considered as the root of
a new finite subtree composed only by paths of the type
$\BbV(\bom{n}^{(\alpha')})$, where $\bom{n}^{(\alpha')}$ is the chromatic matrix
obtained from $\bom{n}^{(\alpha)}$ by changing the root colour from $\alpha$ into 
$\alpha'$. This procedure can be repeated
\textit{ad infinitum} for each leaf and it gives rise to the infinite subtree
$\BbU(\bom{n})$ of $\BbV$ having skeleton in the equivalence class of 
the matrix $\bom{n}$.
 
Notice also that the equivalence relation defined strictly on the space of admissible
chromatic matrices can be extended to more general $d\times d$ matrices of positive entries
differing from admissible ones by terms of order $\cO(1)$.

\subsection{Directional estimates}
Let $\bom{\beta}$ be a $d\times d$ matrix with $\beta_{ij}\geq0$,
$\sum_{i,j=1}^d \beta_{ij}=1$, and 
$\sum_{j=1}^d \beta_{ij}=\sum_{j=1}^d \beta_{ji}$.
For every $n\in\BbN$, denote $\nu_{ij}(n)=\lfloor\beta_{ij} n 
\rfloor$.

\begin{lemm}{}
\label{lem:asympt-K}
Let $\tilde{\bom{\nu}}(n)\in\textsf{ACM}_n$ be an admissible chromatic matrix
in the equivalence class of $\bom{\nu}(n)$. Then, for 
large $n$, we have $K(\bom{\tilde{\nu}}(n))^{1/n}=
\psi(\bom{\beta})(1+o(1))$.
\end{lemm}

\Proof
Notice  that $\lim\limits_{n\r\infty}
\frac{\nu_{ij}(n)}{n}=\beta_{ij}$, hence by the previous remark,
$\lim\limits_{n\r\infty}
\frac{\tilde{\nu}_{ij}(n)}{n}=\beta_{ij}$ also. Now $\tilde{\bom{\nu}}(n)$
is admissible, therefore $K(\tilde{\bom{\nu}}(n))$ is
 obtained by the formula established in lemma \ref{lem:comb-estimates}.
 The proof is
completed by applying Stirling's formula to the factorials appearing in the formula
for $K(\tilde{\bom{\nu}}(n))$.
\eproof

\Rk
The previous matrix $\bom{\beta}$ must be thought as defining
a main direction on the tree. In the sequel,
we shall call \textit{direction} $\bom{\beta}$
a matrix   $\bom{\beta}$ with positive elements
such that $\sum_{i,j=1}^d\beta_{ij}=1$ and verifying
the constraint $\sum_{j=1}^d(\beta_{ij}-\beta_{ji})=0$.

\begin{lemm}{}
For every direction $\bom{\beta}$, there exists
$x_{\bom{\beta}}\in [0,1]$ such that
\[\min\limits_{x\in[0,1]} \phi(x,\bom{\beta})=
\phi(x_{\bom{\beta}},\bom{\beta}).\]
\end{lemm}

\Proof
Since, for every $i,j$, the distribution of the random variables $ \xi_{ij}$ has not 
an atom at 0, it follows that $\BbP( \xi_{ij}>0)=1$. Therefore,
the function $\log m_{ij}(x)=\log\BbE( \xi_{ij}^x)=\log\BbE(\exp(x\log \xi_{ij}))$,
defined on $\BbR$, takes values in $]-\infty, +\infty]$.
Now by H\"older's inequality, it is convex and since
$ \xi_{ij}$ is not a trivial random variable (not almost surely
a constant) it is strictly convex. Hence, the set $\{x\in\BbR: \log m_{ij}(x)<\infty\}$
is an interval of $\BbR$ containing $0$ and $1$ since $\log m_{ij}(0)=0$ and
$\log m_{ij}(1)=\log\BbE(\xi_{ij})<\infty$ by the assumed integrability of $\xi_{ij}$.
Hence $m_{ij}(x)$ is well defined on $[0,1]$ and by standard arguments, infinitely differentiable
in $]0,1[$. 

 Moreover, for $x\in ]0,1[$, $(\log m_{ij})'(x)=\frac{\BbE(\xi_{ij}^x\log\xi_{ij})}{m_{ij}(x)}=
\frac{m'_{ij}(x)}{m_{ij}(x)}$.
Now $\log\phi(x,\bom{\beta})=\log\psi(\bom{\beta})+
\sum_{i,j}
\beta_{ij}\log m_{ij}(x)$; hence
$\frac{\partial}{\partial x} \log\phi(x,\bom{\beta})=
\sum_{i,j}\beta_{ij}\frac{ m'_{ij}(x)}{m_{ij}(x)}$
and either this derivative vanishes for some
$x_{\bom{\beta}}\in ]0,1[$ --- and this
will be a minimum because of convexity ---
or it never vanishes in $]0,1[$ and the
minimum will be attained at one of the borders
of $[0,1]$. Hence the minimum always exists
in $[0,1]$.

Moreover, for non trivial $\xi_{ij}$, 
\textit{i.e.}\ 
not identically equal to a constant for all 
$i,j$, the
minimum is unique.
\eproof

\begin{prop}{}
Let $\bom{\beta}$ be a direction on the tree.
Let $x_{\bom{\beta}}
=\arg\min_{x\in[0,1]}
\phi(x,\bom{\beta})$. Denote by
$\rho(x_{\bom{\beta}})$ the largest eigenvalue
of the matrix $\bom{m}(x_{\bom{\beta}})$ (assumed to be irreducible)
and  $\bom{l}(x_{\bom{\beta}})$ --- 
respectively $\bom{r}(x_{\bom{\beta}})$ --- 
the left --- respectively right ---  eigenvector 
associated with the largest eigenvalue.
On the space of $d\times d$ direction 
matrices, introduce a mapping 
\[S:\bom{\beta}\mapsto \bom{\beta'}=
(\beta'_{ij})_{i,j=1,\ldots,d}\]
by
\[ \beta'_{ij}=b l_i(x_{\bom{\beta}})
m_{ij}(x_{\bom{\beta}}) r_j(x_{\bom{\beta}}).\]
Then
\begin{enumerate}
\item the constant $b$ can be chosen so that 
$\bom{\beta'}$ is also a direction matrix,
\item $\beta'_{ij}=0$ if, and only if, 
$m_{ij}(x_{\bom{\beta}})=0$, and
\item the mapping $S$ has a fixed point.
\end{enumerate}
\end{prop}

\Proof
\begin{enumerate}
\item Remark that $\beta'_{ij}\geq 0$ and since
$\bom{r}$ is a right eigenvector, 
$\sum_j\beta'_{ij}= b l_i(x_{\bom{\beta}})
\rho(x_{\bom{\beta}})r_i(x_{\bom{\beta}})$, we
can chose $b=(\rho(x_{\bom{\beta}})
\bra \bom{l}(x_{\bom{\beta}})\mid 
\bom{r}(x_{\bom{\beta}})\ket)^{-1}$, where 
$\bra\cdot\mid\cdot\ket$ is the scalar product in 
$\BbR^d$, to guarantee $\sum_{ij}\beta'_{ij}=1$.
Similarly, since $\bom{l}$ is
the left eigenvector,
we get 
$\sum_j\beta'_{ji}= b l_i(x_{\bom{\beta}})
\rho(x_{\bom{\beta}})r_i(x_{\bom{\beta}})$, 
fulfilling thus the constraint
$\sum_j(\beta'_{ij}-\beta'_{ji})=0$.
Therefore, $\bom{\beta'}$ is also a direction.
\item
Since $\bom{m}(x)$ is an irreducible matrix, \textit{i.e.} 
there exists $N$ such that $\bom{m}(x)^N$ has 
all its elements strictly positive, the 
eigenvectors corresponding to the
maximal eigenvalue have no zero components.
\item
Let $H$ be the set of direction matrices, viewed as a subset
of the finite dimensional linear space $\BbR^{d^2}$. Then
$H$ is compact and convex and
the map $S:H\r H$ is continuous.
Therefore, by Brouwer's theorem $S$ has (at least)
a fixed point.
\end{enumerate}
\eproof

\begin{prop}{}
If for every direction $\bom{\beta}$, 
$\inf\limits_{x\in[0,1]}\phi(x,\bom{\beta})\leq C$
for some constant $C$, then $\exists x_0\in[0,1]$
such that $\rho(x_0)\leq C$.
\end{prop}
\Proof
Let $\ol{\bom{\beta}}$ be a fixed point of the
map $S$. Then  $\inf\limits_{x\in[0,1]}
\phi(x,\ol{\bom{\beta}})=\phi(y,\ol{\bom{\beta}})
\leq C$.
Compute now $\phi(y,\ol{\bom{\beta}})$ by
using the fact that $\bom{l}$ and $\bom{r}$
are left and right eigenvectors of $\bom{m}$.
\begin{eqnarray*}
\phi(y,\ol{\bom{\beta}})&=&
\left(\frac{b l_1(m_{11}r_1+\ldots+m_{1d}r_d)}
{b l_1 m_{11} r_1} m_{11}\right)^
{\ol{\beta}_{11}}\ldots
\left(\frac{b l_d(m_{d1}r_1+\ldots+m_{dd}r_d)}
{b l_d m_{dd} r_d} m_{dd}\right)^
{\ol{\beta}_{dd}}\\
&=& \rho(y)^{\sum\ol{\beta}_{ij}}
(\frac{r_1}{r_1})^{\ol{\beta}_{11}}
(\frac{r_1}{r_2})^{\ol{\beta}_{12}}
\ldots
(\frac{r_1}{r_d})^{\ol{\beta}_{1d}}
\ldots
(\frac{r_d}{r_1})^{\ol{\beta}_{d1}}
\ldots
(\frac{r_d}{r_d})^{\ol{\beta}_{dd}}
\end{eqnarray*}
and since $\sum_{ij}\ol{\beta}_{ij}=1$
and $\sum_j (\ol{\beta}_{ij}-\ol{\beta}_{ji})
=0$, we have finally that 
\[\phi(y,\ol{\bom{\beta}})=\rho(y)\leq C.\]
\eproof

\begin{coro}{}
If\/ $\inf_{x\in[0,1]}\rho(x)>1$, then there exists a
direction $\ol{\bom{\beta}}$ such that
 $\inf_{x\in[0,1]}\phi(x,\ol{\bom{\beta}})>1$.
\end{coro}
\Proof Immediate from the previous proposition. \eproof

\begin{lemm}{}
\label{lem:continuity}
If\/ $\inf_{x\in[0,1]}\rho(x)>1$, then there exists a
direction $\hat{\bom{\beta}}$ with 
\textit{rational}
coefficients such that
$\inf_{x\in[0,1]}\phi(x,\hat{\bom{\beta}})>1$
\end{lemm}
\Proof
From the previous corollary, there exists a direction
$\ol{\bom{\beta}}$ (with real coefficients)
such that 
$\inf_{x\in[0,1]}\phi(x,\ol{\bom{\beta}})>1$.
Now the function $f(\bom{\beta})=
\inf_{x\in[0,1]}\phi(x,\bom{\beta})$ is continuous.
Hence, for arbitrary directions $\ol{\beta}$ and
$\hat{\beta}$ there exists $C'$ such that
$|f(\ol{\bom{\beta}})
-f(\hat{\bom{\beta}})|\leq C' \|\ol{\bom{\beta}}-\hat{\bom{\beta}}\|$
and since every real can be approximated
arbitrarily well by rationals, the lemma
follows.
\eproof

\subsection{Subtrees, branching and 
Chernoff-Cram\'er bound}
Let $\hat{\bom{\beta}}$ be a direction with
rational coefficients. 
Since all its matrix elements are rational,
there exist 
integers $\gamma$, depending on 
$\hat{\bom{\beta}}$, such that all the products $\gamma\hat{\beta}_{ij}$,
for $i,j=1,\ldots, d$ are integer-valued. 
It is enough to choose the smallest non-zero such integer,
\[\gamma=\inf\{n\in\BbN^*: \lfloor \hat{\beta}_{ij} n
\rfloor = \hat{\beta}_{ij}n=\nu_{ij}\in\BbN,
\forall i,j\in\{1,\ldots,d\}\}.\]

Consider the infinite subtree $\BbU_\gamma(\bom{\nu})$
of $\BbV$ composed    
only from those paths that  
between levels $k\gamma$ and
$(k+1)\gamma$, with $k\in\BbN$, have
prescribed number
$\nu_{ij}=\hat{\beta}_{ij}
\gamma\in\BbN$
of edges of
type $ij$. This subtree --- viewed as a subtree
of $\BbV$ ---  by use of lemma (\ref{lem:asympt-K}), has a
branching ratio (see \cite{Lyo90} for definition
of branching ratio)
$\textrm{br}(\BbU_\gamma(\bom{\nu}))=
\psi(\hat{\bom{\beta}})$, therefore, viewed as a 
tree on its own, it will have branching ratio
$\psi(\hat{\bom{\beta}})^\gamma$.
We give below an extension of a technical result
of \cite{LyoPem} that yields the Chernoff-Cram\'er
bound in our case.
\begin{prop}{}
\label{prop:Chernoff-Cramer}
Suppose $\lambda=\inf_{x\in[0,1]}\rho(x)>1$.
Then there exists a rational direction
$\hat{\bom{\beta}}$, integers $N$ and $\gamma$
in $\BbN^*$, and a real $y\in]0,1]$ such that for all $k\geq N$, 
\[\BbP(\xi[u;v]>y^{\gamma k})>\left(\frac{1}{y 
\psi(\hat{\bom{\beta}})}\right)^{\gamma k},\]
for all $u,v\in\BbU_{\gamma}(\bom{\nu})$ with
$u<v$, $|v|-|u|=k\gamma$, and $\nu_{ij}=
\hat{\beta}_{ij}\gamma=
\lfloor\hat{\beta}_{ij}\gamma\rfloor\in\BbN$.
\end{prop}
\Proof
Since $\lambda>1$, 
by lemma \ref{lem:continuity}
there exists a rational direction 
$\hat{\bom{\beta}}$ such that 
$\inf_{x\in[0,1]} \phi(x,\hat{\bom{\beta}})
>1$. 
Choose $\gamma$
the smallest integer such that all the products
$\hat{\beta}_{ij}\gamma$ are integer-valued 
for all $i,j$.
This choice of $\gamma$ fixes all the integers
$\nu_{ij}=\hat{\beta}_{ij}\gamma =\lfloor \hat{\beta}_{ij}\gamma\rfloor$
for $i,j=1,\ldots, d$. Therefore the subtree
$\BbU_\gamma(\bom{\nu})$ is well defined.
For some $k\in \BbN$, whose value will be fixed later,
let $u,v\in\BbU_{\gamma}(\bom{\nu})$ with $u<v$ and $|v|-|u|=k\gamma$.
The family of the random variables $(\xi_a)_{a\in[u;v]}$ is independent.
Moreover, the product
$\prod\limits_{a\in[u;v]}\xi_a=\xi[u;v]$ can be written
as $\xi[u;v]=\prod\limits_{l=1}^k A_l$, with
$A_l=\prod_{a\in [v|_{|u|+(l-1)\gamma}; 
v|_{|u|+l\gamma}]} \xi_a$.
The family $(A_l)_{l=1,\ldots,k}$ is independent and identically distributed.
Let $\tau(x)=\BbE(A_1^x)$ and $\alpha=\BbE(\log A_1)$ 
and denote $I(x^*)$ the Legendre transform of $\tau(x)$, \textit{i.e.}\
$I(x^*)=\sup_{x\in[0,1]}(xx^*-\log \tau(x))$ for $x^*\in\BbR$. Then the sequence $(S_n)$ with
$S_n=\sum\limits_{l=1}^n \log(A_l)$ satisfies a large
deviation principle with rate function $I$, \textit{i.e.}\
\[\lim_n \frac{1}{n} \log\BbP(S_n > n\theta) 
= - \inf_{x^* > \theta} I(x^*).\]
In other words, for every $\epsilon>0$ we can find an integer $N$ such that
$\frac{1}{n} \log\BbP(S_n > n\theta) 
\geq - \inf_{x^* > \theta} I(x^*)-\epsilon$ for
all $n\geq N$. Choose  any $k\geq N$. Then for this 
sufficiently large $k$, we have
\[\exp(k\epsilon)\BbP(\prod\limits_{l=1}^k A_l \geq (\exp(\theta))^k)
\geq \exp(\inf_{x^*>\theta}I(x^*))^{-k}.\]
This bound is non-trivial if 
$\theta>\alpha$ and then 
$\inf_{x^*>\theta}I(x^*)=I(\theta)$.
Writing, in that case, $\log \hat{y}=\theta$, 
we compute explicitly
$\exp(I(\log\hat{y}))=\frac{\hat{y}}{\hat{y}^{1-x_0} \BbE A_1^{x_0}}$
where $x_0$ is the position where the function $x\log\hat{y}-\log\tau(x)$ attains
its minimum (as a function of $x$). Writing $\hat{y}=\overline{y}^\gamma$, we get
finally
\[
\exp(-k\epsilon)\BbP(\xi[u;v]>\overline{y}^{k\gamma})\geq
\frac{\overline{y}^{k\gamma(1-x_0)}(\BbE A_1^{x_0})^k}{\overline{y}^{k\gamma}}
=\left(\frac{\overline{y}^{1-x_0}\chi(x_0, \hat{\bom{\beta}})}{\overline{y}}\right)^{k\gamma}
=\left(\frac{\overline{y}^{1-x_0}\chi(x_0, \hat{\bom{\beta}})\psi(\hat{\bom{\beta}})}{\overline{y}
\psi(\hat{\bom{\beta}})}\right)^{k\gamma}.
\]
By Fenchel's equality (see \cite{LyoPem}),
\[\sup_{\overline{y}\in]0,1]}\inf_{x\geq0} \overline{y}^{1-x} \chi(x,\hat{\bom{\beta}})
=\inf_{x\in[0,1]}\chi(x,\hat{\bom{\beta}}).\]
Denoting $y$ the value of $\overline{y}\in]0,1]$ attaining the supremum in the
left hand side of Fenchel's equality,
the previous inequality becomes
\[\exp(-k\epsilon)\BbP(\xi[u;v]>y^{k\gamma})\geq
\left(\frac{\inf_{x\in[0,1]}\chi(x,\hat{\bom{\beta}}) \psi(\hat{\bom{\beta}})}
{y \psi(\hat{\bom{\beta}})}\right)^{k\gamma}
=\left(\frac{\inf_{x\in[0,1]}\phi(x,\hat{\bom{\beta}})}
{y \psi(\hat{\bom{\beta}})}\right)^{k\gamma}.\]
Now, $\inf_{x\in[0,1]}\phi(x,\hat{\bom{\beta}})>1$, hence there exists
$\delta>0$ such that $\inf_{x\in[0,1]}\phi(x,\hat{\bom{\beta}})>\exp\delta$.
Choosing then $\epsilon\in]0,\delta\gamma[$, the result follows.
\eproof

\section{Proof of non-ergodicity and transience}
\label{sect:proof}
It will be shown in subsection 
\ref{sect:transience},
that in the case $\lambda>1$ the random walk is transient;
this result is based on some non-probabilistic
estimates coming from the analogy between reversible
Markov chains and electrical networks. We
start however by proving under the same
conditions, in subsection 
\ref{sect:non-ergod}, a seemingly weaker result, 
guaranteeing only non-ergodicity. The reason
is that this method is purely probabilistic
and contains some information that can be
used for the yet unsolved multiplicative
chaos problem for the random walk in random
environment on a coloured tree. 
\subsection{Non-ergodicity for the 
random walk in random environment in the case 
$\lambda>1$}
\label{sect:non-ergod}
Inspired from the result 2.1.7 of \cite{FayMalMen} we prove
the following
\begin{lemm}{}
Let $(X_n)_{n\in\BbN}$ be a real-valued process on 
the filtred space $(\Omega,
\cF,(\cF)_{n\in\BbN},\BbP)$, adapted to the filtration $(\cF_n)_{n\in\BbN}$,
with $X_0=x$, some constant, having uniformly 
bounded increments $Y_n$, \textit{i.e.}  there exists $a>0$, such that
$|Y_n|=|X_n-X_{n-1}|\leq a$ for all $n$. Suppose moreover
that $(X_n)_{n\in\BbN}$ is a strong submartingale \textit{i.e.}\ 
there exists $\epsilon>0$ such that $\BbE(Y_n|\cF_{n-1})\geq
\epsilon$ almost surely for all $n$.
Let $\tau(\delta)=\inf\{n\geq1: X_n<x+n\delta\}$.
Then, there exists some $\delta_1=\delta_1(\epsilon)>0$
such that for all $\delta<\delta_1$, $\BbP(\tau(\delta)=\infty)>0$.
\end{lemm}
\Proof
Write $X_n=x+\sum_{k=1}^nY_k$. Then
\begin{eqnarray*}
\BbP(X_n<x+\delta n) &=& \BbP(-\sum_{k=1}^nY_k>-\delta n)\\
&=& \BbP(\exp(-h\sum_{k=1}^nY_k)>\exp(-h\delta n))\ \ \textrm{for}
\ \ h\geq0\\
&\leq & \exp(h\delta n) \BbE(\exp(-h\sum_{k=1}^nY_k))\\
&=& \exp(h\delta n) 
\BbE(\prod\limits_{k=1}^n\BbE(\exp(-h Y_k)|\cF_{k-1})).
\end{eqnarray*}
Chose $h\leq 1/a$ and use the inequality
$\exp(x)\leq 1+ x+3x^2/2$ valid for $|x|<1$ to bound
the conditional expectation
$\BbE(\exp(-h Y_k)|\cF_{k-1}))\leq 1-h\epsilon+
\frac{3h^2a^2}{2}\equiv \exp(-\delta_2(h))$.
Choosing $h\in]0,\frac{2\epsilon}{3a^2}[$ we get
$\delta_2>0$.
Hence $\BbP(X_n<x+\delta n)\leq \exp(n(h\delta-\delta_2(h)))$ and
choosing $\delta<\delta_1=\delta_2(h)/h$ we get a bound that
is exponentially small for large $n$, that
is $\BbP(X_n<x+\delta n)\leq p_n$ with $p_n=
\exp(-n\delta_3)$ with $\delta_3>0$. Now, by Borel-Cantelli lemma,  
the events
$A_n=\{X_n<x+\delta n\}$ are realised for a finite number
of indices, \textit{i.e.}\  $\forall \gamma \in]0,1[$ we
can find $m=m(\gamma)\in\BbN$ such that $\BbP(\cap_{n=m}^\infty
A_n^c)> \gamma$ and consequently $\BbP(\tau(\delta)=\infty)>0$.
\eproof

We shall now prove that when $\lambda>1$,
not only $\BbE Z_n \r\infty$ but also $Z_n\r\infty$ almost
surely, where 
\[Z_n=\sum_{k=0}^n \sum_{v\in\BbV_k} \xi[v]=
\sum_{k=0}^n Y_k,\]
and this result is enough to prove non-ergodicity for the
random walk in random environment  since it shows that the invariant measure $\pi[v]$ cannot
be normalised. 

\begin{theo}{}
If $\lambda>1$, then  $Z_n\r\infty$ almost surely.
\end{theo}

\Rk This result in conjunction  with the result concerning ergodicity localises
the critical point for the chaos functional equation to $\lambda =1$. 

\noindent\textit{Proof of the theorem:}
From lemma \ref{lem:continuity}, since
$\lambda>1$, there exists a direction with rational
coefficients $\hat{\bom{\beta}}$ such that 
$\inf_{x\in[0,1]}\phi(x, \hat{\bom{\beta}})>1$.
Choose some integer $\gamma$ large enough so that
$\nu_{ij}=\hat{\beta_{ij}}\gamma = 
\lfloor \hat{\beta_{ij}}\gamma \rfloor \in \BbN$ for all
$i$ and $j$ and such that all the integers $\nu_{ij}$ are
so large that Stirling's formula applies.

Consider the infinite subtree $\BbU_\gamma(\bom{\nu})$
having prescribed number of $ij$-type edges between
levels that are at distances $l\gamma$ and $(l+1)\gamma$
from the root. Additionally, since $\lambda>1$, by proposition
\ref{prop:Chernoff-Cramer}, it is possible
to choose $y\in]0,1]$ and
$k\in \BbN$ so that for all $u,v\in\BbU_\gamma(\bom{\nu})$ with
$u<v$ and $|u|=lk\gamma$ and $|v|=(l+1)k\gamma$ for some
$l \in\BbN$, the Chernoff-Cram\'er bound
\[\BbP(\xi[u;v]>y^{k\gamma})>\left(\frac{1}{y\psi(\hat{\bom{\beta}})}
\right)^{k\gamma}\]
applies. We shall construct a minorating process 
$(\tilde{Y}_n)_{n\in\BbN}$ of $(Y_n)_{n\in\BbN}$
such that the corresponding sum process
$\tilde{Z}_n=\sum_{k=0}^n\tilde{Y}_k\r\infty$.

For every $u\in\BbU_\gamma(\bom{\nu})$ define the random set
\[D(u)=\{v\in \BbU_\gamma(\bom{\nu}): v>u; |v|-|u|=k\gamma;
\xi[u;v]>y^{k\gamma}\}.\]
We have
\[y^{k\gamma} \BbE|D(u)|=
y^{k\gamma}\sum_{\stackrel{v\in\BbU_\gamma(\nu)}{v>u; |v|-|u|=\gamma k}}
\BbP(\xi[u;v]>y^{k\gamma})>1.\]

We shall proceed recursively.
Let $\tilde{Y}_0=1$. Define $C_1=B_1=D(\emptyset)$ and
let $B'_1\subseteq B_1$ be such that $1<\sum_{v\in B'_1}y^{|v|}<2$.
Now either $B'_1=\emptyset$ and then we define $\tau=1$ and stop the
process $\tilde{Y}$ or
$B'_1\ne \emptyset$ and then we set $\tau>1$, $B''_1=B_1\setminus B'_1$ and
$\tilde{Y}_1=\sum_{v\in B_1} y^{|v|}$.
On the set $\{\tau>1\}$, we have
$|\tilde{Y}_1-\tilde{Y}_0|
\leq y^{k\gamma} d^{k\gamma} +1<2y^{d\gamma} d^{k\gamma}$. 

Suppose that $\tau>n$ and the sequences $(\tilde{Y}_n), (C_n), (B_n), (B'_n), (B''_n)$
have been constructed up to time $n$. Define $C_{n+1}=\cup_{u\in B'_n} D(u)$,
$B_{n+1}=C_{n+1}\cup B''_n$ and let
$B'_{n+1}\subseteq B_{n+1}$ be an arbitrary subset of $B_{n+1}$ chosen
so that $1<\sum_{v\in B'_{n+1}}y^{|v|}<2$.
Now, either $B'_{n+1}=\emptyset$ and then $\tau=n+1$ and the process
is stopped, or else
$B'_{n+1}\ne \emptyset$ and then $\tau>n+1$, $B''_{n+1}=B_{n+1}\setminus
B'_{n+1}$, $C_{n+1}=\cup_{v\in B'_n} D(v)$,  and 
$\tilde{Y}_{n+1}=\sum_{v\in B_{n+1}} y^{|v|}$.
The increments of the process 
$(\tilde{Y}_{n})$ are bounded
since
\begin{eqnarray*}
|\tilde{Y}_{n+1}-\tilde{Y}_{n}|&=&\left|\sum_{v\in C_{n+1}}y^{|v|}+\sum_{v\in B''_{n}}y^{|v|}
-\sum_{v\in B'_{n}}y^{|v|}-\sum_{v\in B''_{n}}y^{|v|}\right|\\
&=& \left|\sum_{v\in C_{n+1}}y^{|v|}-\sum_{v\in B'_{n}}y^{|v|}\right|\\
&=&\sum_{v\in B'_{n}}y^{|v|}(y^{k\gamma}|D(v)|-1)\\
&\leq& 2 y^{k\gamma} d^{k\gamma}.
\end{eqnarray*}
On the other hand,  the process is a strong submartingale since
the conditional increments with respect to the natural filtration
verify
\[\BbE(\tilde{Y}_{n+1}-\tilde{Y}_{n}|\cF_n)=
\sum_{v\in B'_{n}}y^{|v|}| \BbE(y^{k\gamma}|D(v)|-1|\cF_n)>\epsilon,\]
by virtue of the previous induction step and  of the Chernoff-Cramér bound, 
valid because $\lambda>1$.

Hence, $(\tilde{Y}_l)_{l\in\BbN}$ 
is a strong submartingale
with bounded jumps and by the previous lemma
there exists some index $L$ and some positive $\delta$
such that
$\tilde{Y}_l>1+\delta l$ for all $l\geq L$,
showing that $\tilde{Y}_l\r\infty$ with some
strictly positive probability. A fortiori,
the same conclusion holds for the processes
$(Y_n)$ and $(Z_n)$. Now the event
$\{Z_n\r\infty\}$ is a tail-measurable
event and the random variables 
$(\xi_a)_{a\in\BbA(\BbV)}$
are independent for different generations. Hence
by the 0-1 law, $\BbP(Z_n\r\infty)=1$.
\eproof

\subsection{Transience for the random walk 
in random environment for the case $\lambda>1$}
\label{sect:transience}
To prove transience, we need some not
probabilistic methods based on
electric networks analogy 
(see \cite{Kel} for instance). This is possible
since we are dealing with reversible Markov
chains. 
\begin{theo}{}
If $\lambda>1$, the random walk in random environment
is almost surely transient.
\end{theo}
\Proof
After the long preparatory work on directional estimates, the proof of transience
is essentially reduced to an appropriate extension of the Chernoff-Cram\'er bound.
Recall that a cutset is a finite set, $C$, of vertices not including the root vertex,
 such that any  path from the root vertex to the boundary of the tree 
intersects the set $C$ at exactly one vertex. 
It is enough to show, as it was the case
in \cite{LyoPem}, 
that there exists $w\in\,]0,1[$ such that
for every cutset $C$, there exists some $\epsilon>0$
such that
$\sum_{v\in C}w^{|v|}\xi[v]>\epsilon$.
We only sketch the proof since --- using the
extension of the Chernoff-Cram\'er bound 
established in proposition 
\ref{prop:Chernoff-Cramer} ---
it 
follows the same 
lines as the  proof of item \textit{i)} 
of theorem 1 in \cite{LyoPem}.
Let $\hat{\bom{\beta}}, \gamma, k$, and $y$ be
the parametres determined in  proposition 
\ref{prop:Chernoff-Cramer}. Fix some
sufficiently small $\epsilon>0$ and consider
the random subgraph $\BbW$ of 
$\BbU_\gamma(\bom{\nu})$
whose paths $[u;v]$ with 
$u,v\in \BbU_\gamma(\bom{\nu})$, $u<v$,
and $|v|-|u|=k\gamma$ have been erased
if either $\xi[u;v]<y^{k\gamma}$ holds
or $A_l<\epsilon$ for some $l=1,\ldots,k$ holds,
where $A_l=\xi[v|_{|u|+(l-1)\gamma};
v|_{|u|+l\gamma}]$.
Therefore, by proposition 
\ref{prop:Chernoff-Cramer},
edges remain with probability 
$p>(y\psi(\hat{\bom{\beta}}))^{-k\gamma}$.
Choose $w\in\,](y\psi(\hat{\bom{\beta}})
p^{1/k\gamma})^{-1},1[$.
Since $p\psi(\hat{\bom{\beta}})^{k\gamma}>
1/(wy)^{k\gamma}>1$ by percolation arguments
there is a subtree $\BbW_{k\gamma}$ of $\BbV$
having branching bigger than $1/(wy)^{k\gamma}>1$.
Then we conclude as in \cite{LyoPem}.
\eproof
             
\section{Some results on random walks in
random environments stemming from multiplicative chaos and vice-versa}
\subsection{A multiplicative chaos approach of
the model on a regular tree}
\label{sect:mult-chaos-results}
We use here the results on multiplicative chaos established in a
series of papers \cite{DurLig,Liu97,Liu98} to study the behaviour
of the random walk in random environment on a regular tree
by providing independent proofs of 
those stated in \cite{LyoPem},
stemming solely from results 
on multiplicative chaos.

Let us first remind some basic results on multiplicative chaos. For the problem
of the random walk in random environment on a regular tree,
the multiplicative chaos process reads
\[Y_n[\emptyset]\equiv Y_n=\sum\limits_{v_1,\ldots,v_n=1}^d
\xi_{a(v_1)}\xi_{a(v_1v_2)}\ldots \xi_{a(v_1\ldots v_n)}.\]
Denoting for every $n\in\BbN$ by
$\cF_n=\sigma(\xi_{a(v)}, v\in\cup_{k=0}^n \BbV_k)$,
we have that
\begin{eqnarray*}
\BbE(Y_n|\cF_{n-1}) &=& \sum\limits_{v_1,\ldots,v_{n-1}=1}^d
\xi_{a(v_1)}\xi_{a(v_1v_2)}\ldots \xi_{a(v_1\ldots v_{n-1}n)}
\sum\limits_{v_n=1}^d \BbE\xi_{a(v_1\ldots v_n)}\\
&=& \sum\limits_{v_1,\ldots,v_{n-1}=1}^d
\xi_{a(v_1)}\xi_{a(v_1v_2)}\ldots \xi_{a(v_1\ldots v_{n-1}n)}
\sum\limits_{j=1}^d \BbE\eta_j\\
&=& Y_{n-1} f(1),
\end{eqnarray*}
where we recall that $f(x)=\sum_{i=1}^d \BbE \eta_i^x$ and $g(x)=\log f(x)$.
Now, if $f(1)=1$, the process $(Y_n)$ is a non-negative martingale
converging almost surely to a random variable $Y_\infty\in\BbR^+$.
The interesting and highly non-trivial 
question is whether this convergence holds also in $\cL^1$.

Remark on the other hand that we can write
$Y_n[\emptyset]=\sum_{v\in\BbV_1} \xi_{a(v)} Y'_{n-1}[v]$
where $Y'_{n-1}[v]$ are independent from the $(\xi_{a(v)})_v$. If the
limit in distribution when $n\r\infty$ exists, this gives rise to the functional equation
\begin{equation}
\label{eq:fun-eq-rw}
Y\elaw \sum_{j=1}^d \eta_i Y'_i,
\end{equation}
where $Y$ and $Y'_i$ have the same law for all $i=1,\ldots,d$
and can be chosen independent.  
 If we denote by $\mu$ the common distribution of the random variables
$Y'_i$ and $\hat{\mu}$ its Laplace transform, the functional equation
can be seen as a mapping, $T$, from the space of Laplace transforms of
probability measures into itself, reading
\[T\hat{\mu}(s)=\BbE\left(\prod\limits_{i=1}^d \hat{\mu}(\eta_i s)\right).\]
We use the same symbol, $T$, to denote the induced mapping on the
space of probability measures $\cM_1^+([0,\infty[)$.
Denote $\NTF= \{\mu\in \cM_1^+([0,\infty[): T\mu=\mu\ \ 
\textrm{and}\ \ \mu\ne \delta_0\}$.

There are several known results on the non-trivial fixed points
of $T$ provided some additional moment conditions hold.

\begin{theo} [Durrett and Liggett \cite{DurLig}]
Assume that for some $\delta>1$, $\BbE\eta_i^\delta< \infty$ for all
$i\in\{1,\ldots,d\}$. Then $\NTF\ne \emptyset$ if and only if for some
$\alpha\in ]0,1]$ we have $g(\alpha)=0$ and $g'(\alpha)\leq 0$.
If $g(1)=0$ and $g'(1)<0$, then every $\mu \in \NTF$ 
has $\beta$ moments with $\beta>1$ if and only if
$g(\beta)<0$.
\end{theo}

Liu substantially improved these results in \cite{Liu97,Liu98} by both
weakening the moment condition and by allowing random branching \textit{i.e.}
$d$ becoming a random variable having some known 
joint distribution with the random
variables $(\eta_i)_i$.
We state his result in the special case of interest for us here, namely the case where $d$ is a constant
and where we have assumed that $\BbP(\eta_i>0)=1$ for all $i$.
Notice however that the case of random $d$ allows the treatment of general trees with random branching,
covering --- and even extending --- the ones considered in \cite{LyoPem}. This  study is postponed
to a subsequent paper.

\begin{theo}[Liu \cite{Liu97}]
\label{th:Liu}
Assume that $\BbE[(\sum_{i=1}^d \eta_i) \log^+ (\sum_{i=1}^d \eta_i)]< \infty$,
where $\log^+z=\max(0,\log z)$. 
Then $\NTF\ne\emptyset$ if and only if $g(1)=0$ and
$g'(1)<0$. Moreover, the solutions of the functional equation have finite first moment
and there is a unique such probability measure having first moment equal to 1.
\end{theo}
 
For all subsequent calculations,
it is enough to consider the chaos process
$(Y_n)$ only in the case $f(1)=1$. In the other situations, by renormalising
the random variables we always construct some process satisfying $f(1)=1$.
For this case, 
$(Y_n)$ is a non-negative martingale verifying
$\BbE(Y_n)=1$ for all $n$ and
converging thus almost
surely to a limit $Y_\infty$. By Fatou's lemma, $\BbE(Y_\infty)\leq
\BbE(Y_n)=1$. Thus the limit of the process $(Y_n)$ is always integrable. The question
is whether $\BbE(Y_\infty)=1$. On the other hand, we can seek for solutions of the
functional equation (\ref{eq:fun-eq-rw}).
For such solutions
to be interpreted as limits of the martingale $(Y_n)$, they must be integrable.
By Liu's theorem, we know the conditions of existence of non-trivial integrable
solutions to the functional equation.  The following theorem guarantees
that then the corresponding martingale is uniformly integrable.

\begin{theo}[Kahane and  Peyri\`ere \cite{KahPey}]  
\label{th:KahPey}
With the same notation as above,
the following conditions are equivalent:
\begin{enumerate}
\item $\BbE(Y_\infty)=1$
\item $\BbE(Y_\infty)>0$
\item the functional equation 
(\ref{eq:fun-eq-rw}) 
has an integrable solution 
$Y$ such that $\BbE(Y)=1$.
\end{enumerate}
\end{theo}

We are now able to give the proof of the theorem 
\ref{th-Lyo-Pem} that will be split into three
different \textit{r\'egimes}.

\begin{theo}
Let $\lambda =\inf_{x\in[0,1]} f(x)$. 
If $\lambda<1$, then, almost surely, 
$Z_\infty<\infty$ and the random walk is ergodic.
\end{theo}

\Proof
Since $\lambda<1$, the qualitative behaviour of the function $f$ can only
be in one of the three possibilities depicted in figure \ref{fig:l-lt-1}.

\begin{figure}[h]
\centerline{%
\hbox{%
\psset{unit=4.5mm}
\pspicture(0,0)(39,11) 
\Cartesian(4.5mm,4.5mm) 
\psline{->}(2,2)(12.5,2)
\psline{->}(2,2)(2,11)
\psline[linewidth=1pt,linestyle=dashed]{-}(10,2)(10,2.5)
\psline[linewidth=1pt,linestyle=dashed]{-}(2,4)(12,4)
\psline[linewidth=1pt,linestyle=dashed]{-}(5.6,2)(5.6,4)
\psbezier[linewidth=3pt]{-}(2,9)(5,2)(11,2)(12,3.5)
\rput[rb](1.5,9){$d$}
\rput[rb](1.5,4){$1$}
\rput[tl](5.4,1.5){$\alpha$}
\rput[lt](9.1,1.5){$x_0$}
\rput[tr](12,1.5){$1$}
\psline[linewidth=1pt,linestyle=dashed]{-}(12,2)(12,3.5)
\psline{->}(16,2)(26.5,2)
\psline{->}(16,2)(16,11)
\psline[linewidth=1pt,linestyle=dashed]{-}(23.4,2)(23.4,2.8)
\psline[linewidth=1pt,linestyle=dashed]{-}(16,4)(26,4)
\psbezier[linewidth=3pt]{-}(16,9)(19,2)(25,2)(26,4)
\psline[linewidth=1pt,linestyle=dashed]{-}(19.7,2)(19.7,4)
\rput[lt](19.5,1.5){$\alpha$}
\rput[rb](15.5,9){$d$}
\rput[rb](15.5,4){$1$}
\rput[lt](23.1,1.5){$x_0$}
\rput[tr](26,1.5){$1$}
\psline[linewidth=1pt,linestyle=dashed]{-}(26,2)(26,4)
\psline{->}(28,2)(38.5,2)
\psline{->}(28,2)(28,11)
\psline[linewidth=1pt,linestyle=dashed]{-}(32,2)(32,4)
\psline[linewidth=1pt,linestyle=dashed]{-}(28,4)(38,4)
\psbezier[linewidth=3pt]{-}(28,9)(31,2)(37,2.9)(38,5.5) 
\psline[linewidth=1pt,linestyle=dashed]{-}(34.8,2)(34.8,3.5)
\rput[rb](27.5,9){$d$}
\rput[rb](27.5,4){$1$}
\rput[lt](34.7,1.5){$x_0$}
\rput[lt](32,1.5){$\alpha$}
\rput[tr](38,1.5){$1$}
\psline[linewidth=1pt,linestyle=dashed]{-}(38,2)(38,5.5)
\endpspicture}}
\caption{
\label{fig:l-lt-1}
\textsf{The generic behaviour of the function $f$ in the subcritical case $\lambda<1$.}}
\end{figure}

In all these cases, there is an $\alpha\in]0,1[$ such that $g(\alpha)=0$ and $g'(\alpha)<0$.
Introduce now the renormalised
random variables $\tilde{\eta}_i=\eta^\alpha_i$ for this value of $\alpha\in]0,1[$.
We denote with a tilde all the objects relative to the new random variables $\tilde{\eta_i}$,
like $\tilde{\xi}[v]$, $\tilde{Y}_n$, $\tilde{f}$, or $\tilde{g}$,
defined in the same way as the corresponding objects without tilde for the
random variables $\eta_i$.
We compute for instance $\tilde{f}(x)=\sum_{i=1}^d \tilde{\eta}_i^x=f(\alpha x)$
and $\tilde{g}(x)=g(\alpha x)$.

Introduce now the variables $\hat{\eta}_i=\frac{\tilde{\eta}_i^t}{\tilde{f}(t)}$
for some value of $t$ that will be determined later.
All the quantities relative to the new variables $\hat{\eta}_i$ will now be distinguished
by the caret symbol $\mbox{}\ \ \hat{}\ \ $. We have
\[\hat{g}(y)=g(\alpha t y)- yg(\alpha t)\]
and
\[ \hat{g}'(y)=\alpha t g'(\alpha t y)- g(\alpha t),\]
establishing thus that $\hat{g}(1)=0$ for all possible choices of $t$.

Let $\beta_0$ be the abscissa of the point of the graph of $g$ that
lies also on the straight line from the origin tangent to the graph
(see figure
\ref{fig:l-lt-1-tangent}).


\begin{figure}[h]
\centerline{%
\hbox{%
\psset{unit=4.5mm}
\pspicture(0,0)(14,12)
\Cartesian(4.5mm,4.5mm)
\psline{->}(2,5)(12,5)
\psline{->}(2,2)(2,11)
\psline[linewidth=1pt,linestyle=dashed]{-}(9.8,2.5)(9.8,5)
\psline[linewidth=1pt,linestyle=dashed]{-}(7.8,2.9)(7.8,5)
\psbezier[linewidth=3pt]{-}(2,9)(5,2)(11,2)(12,3.5)
\psline{-}(2,5)(12,1.5)
\rput[rb](1.5,9){$\log d$}
\rput[rb](1.5,5){$0$} 
\rput[bl](5.2,5.5){$\alpha$}
\rput[bl](9.5,5.5){$x_0$}
\rput[bl](7.5,5.5){$\beta_0$}
\endpspicture}}
\caption{
\label{fig:l-lt-1-tangent}
\textsf{
The graph of the fucntion $g(x)=\log f(x)$ in the case $\lambda<1$.
The straight line from the origin tangent to the graph of $g$ touches
the graph at a point with abscissa $\beta_0\in\,]\alpha,x_0[$.
}}
\end{figure}

Choose now $t=\beta/\alpha$ for some  $\beta\in ]\alpha,\beta_0[$.
For that choice of $t$, we have $g'(1)=\beta g'(\beta)-g(\beta)<0$ and
$\tilde{f}(t)=f(\beta)<f(\alpha)<1$.

Since now $\hat{g}(1)=0$ and $\hat{g}'(1)<0$
the process $(\hat{Y}_n)_n$ 
with
\[\hat{Y}_n=
\sum\limits_{v\in\BbV_n}\hat{\xi}[v]=
\frac{1}{f(\beta)^n}
\sum\limits_{v\in\BbV_n}\hat{\xi}[v],\]
converges in $\cL^1$ to a random
variable $\hat{Y}_\infty\in\BbR^+$ by  theorem 
\ref{th:Liu}.
Since $f(\beta)<1$ and 
$\hat{Y}_\infty<\infty$ almost surely, then
\[\sum_n \sum_{v\in\BbV_n} \xi^\beta[v]=\sum_n f(\beta)^n \hat{Y}_n <\infty \ \  
\textrm{almost surely.}\]
Hence we conclude as in lemma  \ref{lem:ergod} 
that $Z_\infty<\infty$ almost surely.
Therefore the invariant measure is normalisable and
the random walk is ergodic.
\eproof

\begin{theo}
If $\lambda=\inf_{x\in[0,1]}f(x)=1$ and $\sum_{i=1}^d \BbE(\eta_i \log \eta_i)<0$
then, almost surely,  $0<Y_\infty<\infty$,
$Z_\infty=\infty$, and the random walk is null-recurrent.
\end{theo}
\Proof
Denote by $x_0$ the unique point of the interval $]0,1]$ attaining the infimum, \textit{i.e.}
$f(x_0)=\lambda=1$. By the strict convexity of the function $g(x)=\log f(x)$, we
can have only three possibilities, depicted generically for the function $f$  on the figure \ref{fig:l-eq-1} below;
similar figures can be drawn for $g=\log f$.


\begin{figure}[h]
\centerline{%
\hbox{%
\psset{unit=4.5mm}
\pspicture(0,0)(39,11)
\Cartesian(4.5mm,4.5mm)
\psline{->}(2,2)(12,2) 
\psline{->}(2,2)(2,11)
\psline[linewidth=1pt,linestyle=dashed]{-}(9.4,2)(9.4,4)
\psline[linewidth=1pt,linestyle=dashed]{-}(2,4)(12,4)
\psbezier[linewidth=3pt]{-}(2,9)(7.5,2)(11,4)(12,5)
\rput[rb](1.5,9){$d$}
\rput[rb](1.5,4){$1$}
\rput[lt](9.1,1.5){$x_0=1$}
\psline{->}(16,2)(26,2)   
\psline{->}(16,2)(16,11)
\psline[linewidth=1pt,linestyle=dashed]{-}(23.4,2)(23.4,4)
\psline[linewidth=1pt,linestyle=dashed]{-}(16,4)(26,4)
\psbezier[linewidth=3pt]{-}(16,9)(20,5.2)(22.5,4.2)(26,3.2)
\rput[rb](15.5,9){$d$}
\rput[rb](15.5,4){$1$}
\rput[lt](23.1,1.5){$x_0=1$}
\psline{->}(28,2)(38,2)
\psline{->}(28,2)(28,11)
\psline[linewidth=1pt,linestyle=dashed]{-}(32,2)(32,4)
\psline[linewidth=1pt,linestyle=dashed]{-}(28,4)(38,4)
\psbezier[linewidth=3pt]{-}(28,9)(31,2.9)(31.5,2.9)(38,7.5)
\rput[rb](27.5,9){$d$} 
\rput[rb](27.5,4){$1$} 
\rput[lt](32,1.5){$x_0<1$}
\endpspicture}}
\caption{
\label{fig:l-eq-1}   
\textsf{The generic behaviour of the function $f$ in the critical case $\lambda=1$.
In a) $f(1)=1$ and $f'(1)=0$, in b) $f(1)=1$ and $f'(1)<0$, and in c)
$f(x_0)=1$ but $f(1)>1$ and $f'(1)>0$.}}
\end{figure}

Cases a) and c) are excluded by 
the condition $\sum_{i=1}^d 
\BbE(\eta_i \log \eta_i)<0$
so that it is enough to consider the case b)  where $g(1)=0$ and $g'(1)<0$.
Hence by  theorem \ref{th:Liu}, 
the functional equation (\ref{eq:fun-eq-rw})
 has non-trivial integrable
solutions
and by theorem \ref{th:KahPey}
the process $(Y_n)_n$ is an 
uniformly integrable martingale converging
in $\cL^1$ to the random variable $Y_\infty$, \textit{i.e.}
$\lim_{n\r\infty}\BbE|Y_n-Y_\infty|=0$  
and subsequently the Ces\`aro mean,
$\frac{1}{n}\sum_{k=0}^{n-1}(Y_k-Y_\infty)\r 0$.
Therefore, $Z_n/n\r Y_\infty$ almost surely and 
$Y_\infty\in \BbR^+\setminus\{0\}$
almost surely because the martingale is uniformly integrable,
establishing thus that $Z_\infty=\infty$ almost surely. This guarantees that the walk is not
ergodic.
On the other hand, since 
$Y_n\r Y_\infty\in \BbR^+\setminus\{0\}$ 
almost surely, using the Nash-Williams criterion based on the  electric circuit analogy
(see corollary 4.3 of \cite{Lyo90} for instance) 
the
walk  cannot be  transient. Hence the walk is almost surely null-recurrent.
\eproof

Finally we treat the supercritical case $\lambda>1$.

\begin{theo}
Let $\lambda =\inf_{x\in[0,1]} f(x)$.
If $\lambda>1$ then almost surely $Y_\infty=\infty$,
$Z_\infty=\infty$, and the random walk is
transient.
\end{theo}
\Proof
Let $x_0\in[0,1]$ be such that  
$f(x_0)=\inf_{x\in[0,1]}f(x)=\lambda>1$.
The qualitative behaviour of the function $f$ has two possibilities
depicted in figure (\ref{fig:l-gt-1}) below.

\begin{figure}[h]
\centerline{%
\hbox{%
\psset{unit=4.5mm}
\pspicture(0,0)(27,11)
\Cartesian(4.5mm,4.5mm)
\psline{->}(2,2)(12,2)
\psline{->}(2,2)(2,11)
\psline[linewidth=1pt,linestyle=dashed]{-}(2,4)(12,4)
\psline[linewidth=1pt,linestyle=dashed]{-}(12,2)(12,5)
\psbezier[linewidth=3pt]{-}(2,9)(5,6)(11,5)(12,5)
\rput[rb](1.5,9){$d$}
\rput[rb](1.5,4){$1$}
\rput[rt](12,1.5){$x_0=1$}
\psline{->}(16,2)(26,2)
\psline{->}(16,2)(16,11)
\psline[linewidth=1pt,linestyle=dashed]{-}(22.5,2)(22.5,5.1)
\psline[linewidth=1pt,linestyle=dashed]{-}(16,4)(26,4)
\psbezier[linewidth=3pt]{-}(16,9)(19,4)(25,4)(26,7)
\rput[rb](15.5,9){$d$}
\rput[rb](15.5,4){$1$}
\rput[lt](22,1.5){$x_0<1$}
\endpspicture}}
\caption{
\label{fig:l-gt-1}
\textsf{
The generic behaviour of the function $f$ in the supercritical case $\lambda>1$.}}
\end{figure}

We have omitted the totally trivial case where the infimum is attained at $x_0=0$, in
which case $g(0)=\log d$ and $g'(0)\geq0$ since then necessarily $\sum_{i=1}^d\eta_i\geq 1$
and consequently $Y_n\r\infty$ exponentially fast.
In case a) where $x_0=1$, 
introduce the random variables $\tilde{\eta}_i=\frac{\eta_i}{f(1)}$
with corresponding function $\tilde{g}(x)=g(x)-xg(1)$. We have
$\tilde{g}(1)=0$ and $\tilde{g}'(1)=g'(1)-g(1)<0$ since $g'(1)\leq 0$
and $g(1)>0$. Thus by theorem \ref{th:Liu}
we have again the uniform integrability of the 
martingale $(\tilde{Y}_n)_n$ reading
\[\tilde{Y}_n =\frac{1}{\lambda^n} \sum\limits_{v\in\BbV_n} \xi[v] =\frac{Y_n}{\lambda^n},\]
which converges in $\cL^1$ to a random variable with $0<\tilde{Y}_\infty<\infty$.
Since $\lambda>1$, necessarily $Y_n\r\infty$ and a fortiori $Z_\infty=\infty$ almost
surely.

In case b) where $x_0<1$, we consider the random variables 
$\tilde{\eta}_i=\frac{\eta_i^{x_0}}{f(x_0)}$ with $\tilde{g}(x)=g(xx_0)-xg(x_0)$
and the random variables 
$\hat{\eta_i}=\frac{\tilde{\eta}_i^t}{\tilde{f}(t)}$ with
$\hat{g}(y)=g(tyx_0)-y(g(tx_0)$ for some $t$ that will be determined later.
Obviously $\hat{g}(1)=0$ for all choices of $t$. Compute
$\hat{g}'(1)=tx_0g'(tx_0)-g(tx_0)$ and define $\beta_0$ as the unique
point in $]0,1]$ that is abscissa  of the point where a straight line
emanating from the origin is tangent to the graph of $g$ (see figure
\ref{fig:l-gt-1-tangent}).

 
\begin{figure}[h]
\centerline{%
\hbox{%
\psset{unit=4.5mm}
\pspicture(0,0)(14,12)
\Cartesian(4.5mm,4.5mm)
\psline{->}(2,2)(12,2)
\psline{->}(2,2)(2,11)
\psline[linewidth=1pt,linestyle=dashed]{-}(9.5,2)(9.5,4.2)
\psline[linewidth=1pt,linestyle=dashed]{-}(7.8,2)(7.8,3.8)
\psbezier[linewidth=3pt]{-}(2,9)(3,2)(11,3)(12,6)
\psline{-}(2,2)(12,4.5)
\rput[rb](1.5,9){$\log d$}
\rput[rb](1.5,2){$0$}
\rput[tl](9.5,1.5){$\beta_0$}
\rput[tl](7.5,1.5){$x_0$}
\endpspicture}}
\caption{
\label{fig:l-gt-1-tangent}
\textsf{The graph of the fucntion $g(x)=\log f(x)$ in the case $\lambda>1$.
The straight line from the origin tangent to the graph of $g$ touches
the graph at a point with abscissa $\beta_0>x_0$.
}}
\end{figure}

 Necessarily,
$\beta_0>x_0$.
Choose an arbitrary $\beta\in ]0,\beta_0[$ 
and let $t=\beta/x_0$.
For that choice of $t$, $\hat{g}'(1)=\beta g'(\beta)-g(\beta)$ and since
$\beta<\beta_0$, we have that $\hat{g}'(1)<0$. Observe also that
$\hat{\eta_i}=\frac{\eta_i^\beta}{f(\beta)}$.
Hence the process $(\hat{Y}_n)_n$ with
\[\hat{Y}_n=\sum\limits_{v\in\BbV_n}\frac{\xi^\beta[v]}{f(\beta)^n}\]
is a uniformly integrable martingale converging in $\cL^1$ to an almost
surely positive limit.

Denote $Y_n^{(\beta)}=\sum\limits_{v\in\BbV_n}\xi^\beta[v]$. Since
$\frac{Y_n^{(\beta)}}{f(\beta)^n}\r \hat{Y}_\infty$ in $\cL^1$ with
$0<\hat{Y}_\infty<\infty$ for all $\beta<\beta_0$,
it follows that
$\frac{1}{n}\log Y_n^{(\beta)} \r g(\beta)$ for all $\beta<\beta_0$.

Now by  a standard argument used in statistical mechanics
(subadditivity and law of large numbers)
for all $\beta\in[0,1]$, $\liminf\frac{1}{n}\log Y_n^{(\beta)}=\sigma(\beta)$ 
where $\sigma(\beta)$ is a convex function.
For $\beta<\beta_0$ the sequence has a limit 
(hence it coincides with its \textit{limes
infimum}) 
so that $\sigma(\beta)=g(\beta)$ for $\beta<\beta_0$. 
By convexity, the graph of  $\sigma(\beta)$ for $\beta>\beta_0$ is bounded
from below by the tangent of $g(\beta)$ at $\beta_0$ (this argument is used
in the context of statistical mechanics in general \cite{Rue} and that
of disordered systems in particular in \cite{ColKou}).
Since $g'(\beta_0)>0$ and $g(\beta_0)>\lambda$,
we have that
$\sigma(\beta)\geq \lambda +(\beta-\beta_0) g'(\beta_0)>0$ for all $\beta>\beta_0$.
Hence $Y_n^{(\beta)}\r\infty$ for all values $\beta>\beta_0$ and in particular
for $\beta=1$ so that $Y_n$ tends to $\infty$ exponentially fast and a fortiori $Z_\infty=\infty$.
Since the process $Y_n$ diverges exponentially fast, using again Nash-Williams criterion,
for every $\epsilon>0$ we
can find some $w\in]0,1[$ such that for every cutset $C$,
$\sum_{v\in C}w^{|v|}\xi[v]>\epsilon$ and thus the walk is transient.
\eproof

\subsection{Some open problems on multiplicative
chaos}

The last subsection demonstrated the close relationship between results on
multiplicative chaos and reversible Markov chains. In particular, the most
difficult part for the Markov chain problem, namely the
critical case $\lambda=1$ becomes an immediate consequence of the theorem
on the existence of non trivial solutions of the functional equation and
the uniform integrability of the corresponding martingale, once the
conditions for the existence of non trivial solutions are known.
This analogy can even be extended on more general settings to include
the case of random trees and of general distributions for the environment
that correspond  to situations much more general than the one considered
in \cite{LyoPem}. Actually, what plays an important \textit{r\^ole} is
the theorem (1) of \cite{DurLig} but this theorem is properly generalised
by Liu \cite{Liu98}
to include random number of variables $d$. Therefore, the treatment
of random walks in general random environment on random trees 
becomes accessible by virtue of the results of Liu  on multiplicative chaos.

We got conditions under which the chaos processes $Y_n$ and $Z_n$ tend to $\infty$ or remain
finite according
to the values of the parametre $\lambda$. The precise study of this
classification gives rise 
to a multiplicative chaos functional equation of the
type
\[Y^{(\alpha)}\ \elaw\ \sum_{ \beta}\eta_{\alpha\beta} Y^{'(\alpha \beta)}\]
for which the conditions of  existence of non trivial
solutions  are not known. 
In view of the results on the random walk problems it is expected
that the classifying parametre in this problem is  the largest eigenvalue of the matrix of moments
$\bom{m}(x)$. This problem is actually under investigation. The
above mentioned intuition is confirmed by some preliminary results, by the
partial results of \cite{BenNas} and by physical intuition. As a matter of
fact the random walk in a random environment can also be viewed
as a physical system of spins in a quenched disorder. In the random string
problem the quenching is quite stringent so that
the Lyapunov's exponent appear. On the contrary, the random walk
in random environment on the coloured tree
behaves very much like a self-averaging problem. 

\vskip5mm
\noindent
\textbf{Acknowledgments:}
M.M.\ wishes to thank Yasunari Higuchi for useful discussions and D.P.\ \textsc{fapesp} for financial support
and the Instituto de Matem\'atica e Estat\'{\i}stica, University of S\~ao Paulo 
for warm hospitality during the period when this paper was prepared.
\bibliographystyle{plain}
\scriptsize{
\bibliography{rwre}
 }            

\end{document}